\documentclass{article}
\usepackage{amssymb}
\usepackage{amsmath}
\usepackage{amsthm}
\usepackage[english]{babel}
\usepackage{authblk}
\usepackage{graphicx}
\newtheorem{theorem}{Theorem}
\newtheorem{lemma}{Lemma}
\newtheorem{prop}{Proposition}
\newtheorem{definition}{Definition}
\newtheorem{remark}{Remark}
\usepackage[margin=2.6cm]{geometry}
\usepackage{soul}

\usepackage{xcolor}
\newcommand{\cred}[1]{\textcolor{red}{#1}}

\newcommand{\iu}{\hat \imath}
\newcommand{\mi}[1]{\mathbf{#1}} 
\newcommand{\vet}[1]{\mathbf{#1}}   
\DeclareMathOperator{\diag}{diag}
\newenvironment{prova}{\noindent {\bf Proof. }}{\hfill $\bullet$ \\}

\begin{document}

\title{Asymptotic spectral properties and preconditioning of an approximated nonlocal Helmholtz equation with fractional Laplacian and variable coefficient wave number $\mu$ }
\author[1]{Andrea Adriani\thanks{andrea.adriani@uninsubria.it}}
\author[2]{Rosita L. Sormani\thanks{rl.sormani@uninsubria.it}}
\author[3]{Cristina Tablino-Possio\thanks{cristina.tablinopossio@unimib.it}}
\author[4,6]{Rolf Krause \thanks{rolf.krause@usi.ch}}
\author[1,5]{Stefano Serra-Capizzano\thanks{stefano.serrac@uninsubria.it}}

\affil[1]{Department of Science and High Technology, University of Insubria, Via Valleggio 11, 22100 Como, Italy,}
\affil[2]{Department of Theoretical and Applied Sciences, University of Insubria, Via Dunant 3, 21100 Varese, Italy}
\affil[3]{Department of Mathematics and Applications, University of Milano - Bicocca, via Cozzi, 53,  20125  Milano, Italy,}
\affil[4]{ Center for Computational Medicine in Cardiology, University of Italian Switzerland, Via Giuseppe Buffi, 13, 6900 Lugano, Switzerland}
\affil[5]{Division of Scientific Computing, Department of Information Technology, Uppsala University, Lägerhyddsv 2, hus 2, SE-751 05, Uppsala, Sweden}
\affil[6]{UniDistance Suisse, Schinerstrasse 18, 3900 Brig, Switzerland}

\date{}

\maketitle
\begin{abstract}
The current study investigates the asymptotic spectral properties of a finite difference approximation of nonlocal Helmholtz equations with a fractional Laplacian and a variable coefficient wave number $\mu$, as it occurs when considering a wave propagation in complex media, characterized by nonlocal interactions and spatially varying wave speeds. More specifically, by using tools from Toeplitz and generalized locally Toeplitz theory, the present research delves into the spectral analysis of nonpreconditioned and preconditioned matrix sequences, with the main novelty regarding a complete picture of the case where $\mu=\mu(x,y)$ is nonconstant. We report numerical evidence supporting the theoretical findings. Finally, open problems and potential extensions in various directions are presented and briefly discussed.
\end{abstract}

\section{Introduction}

In the present work, we consider the following two-dimensional nonlocal Helmholtz equation
\begin{equation}\label{a1}
    \begin{cases}
        (-\Delta)^{\alpha/2}u(x,y)+\mu(x,y)u(x,y)=v(x,y),
        & (x,y)\in \Omega \subset \mathbb{R}^{2},\quad\alpha\in (1,2), \\
        u(x,y)=0, & (x,y)\in\Omega^{c},
    \end{cases}
\end{equation}
with a given variable-coefficient complex-valued wave number $\mu(x,y)$ and source term $v(x,y)$. The fractional Laplacian operator $(-\Delta)^{\alpha/2}(\cdot)$ has the explicit expression
\begin{equation*}
    (-\Delta)^{\alpha/2}\big(u(x,y)\big) =
    c_{\alpha}\text{P.V}.\int_{\mathbb{R}^{2}}\frac{u(x,y)-u(\tilde{x},\tilde{y})}{[(x-\tilde{x})^{2}+(y-\tilde{y})^{2}]^{\frac{2+\alpha}{2}}}d\Tilde{x}d\Tilde{y},
    \qquad c_{\alpha} = \frac{2^{\alpha}\Gamma(\frac{\alpha+2}{2})}{\pi|\Gamma(\frac{-\alpha}{2})|},
\end{equation*}
{in which P.V. stands for the Cauchy principal value} and $\Gamma$ is the Gamma function.
For simplicity, here the domain $\Omega$ is the square $[0,1]^2$ and $\Omega^{c}$ is the complement of $\Omega$ in $\mathbb{R}^{2}$. However, we discuss at the end of the present paper that our analysis can be generalized to any rectangle and also to non-Cartesian domains.

The motivation for this work arises from the fact that, in order to accurately and finely describe certain physical processes such as complex and anomalous diffusion phenomena, fractional operators and fractional differential equations (FDEs) are required. Due to the richness of these models, a great amount of research has been devoted in the last three decades to their connection with anomalous processes. For instance, a detailed study of the analytic aspects can be found in \cite{model-1,model-2,model-3,model-4} and references therein. A particularly versatile operator is the fractional Laplacian: its wide range of applications include the modeling of long-range interactions and of anomalous diffusion \cite{9,10}.

On the downside, as often happens in applied mathematics, analytical solutions are not usually available in explicit form or, when they are, they are computationally inefficient. As a consequence, many numerical methods for approximating the continuous problems have been designed or adapted to the specific setting, such as finite difference-quadrature \cite{13,20,14,15}, spectral methods \cite{16} or ad hoc finite element methods \cite{17,19,18}. Such a variety of research papers is not redundant, since there exist several equivalent definitions of the fractional Laplacian \cite{11,12} that can be treated in different ways, depending also on the desired precision. The computational bottleneck comes once the approximated equations are obtained, usually in the form of a large and structured linear system: at this point, fast methods for their solution are strongly needed.

For this purpose, a circulant preconditioner was developed by Hao, Zhang, and Du \cite{20} for the two-level Toeplitz matrix structure arising from \eqref{a1}. Subsequently, two varieties of $\tau$ preconditioners were proposed by Li, Chen, H.W. Sun, and T. Sun \cite{tau-prec}, who also reported a localization analysis for the spectra of the associated preconditioned matrix sequences. Our contribution lies in the latter direction, motivated by the well-known fact that $\tau$ preconditioning leads to superior performances than other matrix algebra approaches, such as $\omega$-circulants with $|\omega|=1$, when real symmetric Toeplitz structures are involved. We refer the reader to \cite{tau-essential,tau-theory1,tau-theory2} for theoretical studies in a general Toeplitz setting and to \cite{BEV} for a recent interesting contribution in the fractional setting.

The aim of the current paper is to complete the analysis in \cite{tau-prec} by giving a global spectral and singular value analysis of the original coefficient matrix sequences and of the preconditioned ones. We rely on the generalized locally Toeplitz (GLT) apparatus \cite{GLT-BookI,GLT-BookII} and on results concerning the extreme eigenvalues of Toeplitz matrix sequences \cite{BG-extr,S-extr2,S-extr1}. 
More in detail, in \cite{tau-prec} only the case of a constant $\mu$ is theoretically considered, while the specific case of a real-valued variable coefficient $\mu(x,y)=\cos x \cos y$
was considered in the numerical tests: here we provide an exhaustive theoretical analysis of the spectral features of the considered matix-sequences, both in a preconditioned and nonpreconditioned setting, in the general case where $\mu(x,y)$ is nonconstant and complex-valued. Since the presence of  a variable-coefficient complex-valued wave number spoils both the Toeplitzness and the real symmetry of the matrices, specific tools are also needed (see \cite{GLT-BookI,GLT-BookII,orth-distribution2} and references therein).

The paper is organized in the following manner. Section \ref{sec:sp-tools} is an introduction to the theoretical apparatus and a summary of the necessary tools. Section \ref{sec:appr} briefly describes the approximation of the considered equation and emphasises the structural features of the resulting matrices, including a detailed study of the singular value and spectral properties via the GLT tools. In Section \ref{sec:sp-analysis}, the singular value and eigenvalue distribution of the preconditioned matrix sequences is obtained. Finally, in Section \ref{sec:num-exp} the theoretical results are supported by several numerical experiments and in Section \ref{sec:final} conclusions, observations, and open problems are presented.

\section{Spectral tools}\label{sec:sp-tools}

In this introductory section, we present the tools needed to perform the spectral analysis of the involved matrices, derived from the theory of multilevel block GLT matrix-sequences. The non-block setting, corresponding to scalar-valued symbols, is described in detail in books \cite{GLT-BookI,GLT-BookII}, while the block case, corresponding to matrix-valued symbols, is explored in papers \cite{GLT-block1D,GLT-blockdD}. In our specific context, the size of the blocks is $r=1$ and the dimensionality of the domain of the problem is $d=2$, therefore throughout this work we remain in the setting of 2-level non-block GLT sequences. However, the spectral tools are presented in the maximal generality required by the potential extensions.


The preliminary theory is organized as follows. Subsection \ref{ssez:notation} fixes notations and terminology. Subsection \ref{ssez:matrix-seq} introduces the definitions of distribution and clustering for a generic matrix sequence, while Subsection \ref{ssez:acs} delves into the related notion of approximating class of sequences. Subsection \ref{ssez:Toeplitz-etc} is devoted to the matrix structures that make up the multilevel block GLT $*$-algebra, which is treated in Subsection \ref{ssez:GLT} from the point of view of the related operative features.

\subsection{Notation and terminology}\label{ssez:notation}

\paragraph{Matrices and matrix sequences.}
Given a square matrix $A\in\mathbb{C}^{m\times m}$, we indicate with $A^*$ its conjugate transpose and with $A^\dag$ the Moore--Penrose pseudoinverse of $A$. Recall that $A^\dag=A^{-1}$ whenever $A$ is invertible. The singular values and eigenvalues of $A$ are denoted respectively by $\sigma_1(A),\ldots,\sigma_m(A)$ and $\lambda_1(A),\ldots,\lambda_m(A)$.

Regarding matrix norms, $\|\cdot\|$ refers to the spectral norm and for $1\le p\le\infty$ the notation $\|\cdot\|_p$ stands for the Schatten $p$-norm, defined as the $p$-norm of the vector of the singular values. Note that the Schatten $\infty$-norm, which is equal to the largest singular value, coincides with the spectral norm $\|\cdot\|$; the Schatten 1-norm, since it is the sum of the singular values, is often referred to as the trace-norm; and the Schatten 2-norm coincides with the Frobenius norm. Schatten $p$-norms, as important special cases of unitarily invariant norms, are treated in detail in a wonderful book by Bhatia \cite{Bhatia-book}.

Finally, the expression \textit{matrix sequence} identifies any sequence of the form $\{A_n\}_n$, where $A_n$ is a square matrix of size $d_n$ with $d_n$ strictly increasing, so that $d_n\to\infty$ as $n\to\infty$. A $r$-\textit{block matrix sequence}, or simply a matrix sequence if $r$ can be deduced from context, is a special $\{A_n\}_n$ in which the size of $A_n$ is $d_n=r\phi_n$, with $r\ge 1\in\mathbb{N}$ fixed and $\phi_n\in\mathbb{N}$ strictly increasing.

\paragraph{Multi-index notation.}

To effectively deal with multilevel structures it is necessary to use multi-indexes, which are vectors of the form $\mi{i}=(i_1,\ldots,i_d)\in \mathbb{Z}^d$. The related notation is listed below.
	\begin{itemize}
		\item \textbf{0},\textbf{1},\textbf{2},... are vectors of all zeroes, ones, twos,...
		\item $\mi{h}\leq\mi{k}$ means that $h_r\leq k_r$ for all $r=1,\ldots,d$. More in general, relations between multi-indexes are evaluated componentwise.
  		\item Operations between multi-indexes, such as addition, subtraction, multiplication and division, are also performed componentwise.
		\item The multi-index interval $[\mi{h},\mi{k}]$ is the set $\{\mi{j}\in\mathbb{Z}^d:\mi{h}\leq\mi{j}\leq\mi{k}\}$. We always assume that the elements in an interval $[\mi{h},\mi{k}]$ are ordered in the standard lexicographic manner:
		$$\left[\cdots\Big[\big[(j_1,\ldots,j_d)\big]_{j_d=h_d\ldots,k_d}\Big]_{j_{d-1}=h_{d-1},\ldots,k_{d-1}}\cdots\right]_{j_1=h_1,\ldots,k_1}.$$
		\item $\mi{j}=\mi{h},\ldots,\mi{k}$ means that $\mi{j}$ varies from $\mi{h}$ to $\mi{k}$, always following the lexicographic ordering.
        \item $\mi{m}\to \infty$ means that $\min (\mi{m})=\min_{j=1,\ldots,d}m_j\to\infty$.
		\item The product of all the components of $\mi{m}$ is denoted as $\nu(\mi{m}) := \prod_{j=1}^{d}m_j$.
	\end{itemize}
A multilevel matrix sequence is a matrix sequence $\{A_{\mi{n}}\}_n$ such that $n$ varies in some infinite subset of $\mathbb{N}$, $\mi{n} = \mi{n}(n)$ is a multi-index in $\mathbb{N}^d$ depending on $n$, and $\mi{n}\to \infty$ when $n\to \infty$. This is typical of many approximations of differential operators in $d$ dimensions.


\paragraph{Measurability.}
All the terminology from measure theory, such as ``measurable set'', ``measurable function'', ``a.e.'', etc., refers to the Lebesgue measure in $\mathbb R^t$, denoted with $\mu_t$. A matrix-valued function $f:D\subseteq\mathbb R^t\to\mathbb{C}^{r\times r}$ is said to be measurable (resp., continuous, Riemann-integrable, in $L^p(D)$, etc.) if all its components $f_{\alpha\beta}:D\to\mathbb{C},\ \alpha,\beta=1,\ldots,r$, are measurable (resp., continuous, Riemann-integrable, in $L^p(D)$, etc.). If $f_m,f:D\subseteq\mathbb R^t\to\mathbb C^{r\times r}$ are measurable, we say that $f_m$ converges to $f$ in measure (resp.,
a.e., in $L^p(D)$, etc.) if $(f_m)_{\alpha\beta}$ converges to
$f_{\alpha\beta}$ in measure (resp., a.e., in $L^p(D)$, etc.) for all $\alpha,\beta=1,\ldots,r$.

\subsection{Distribution and clustering}\label{ssez:matrix-seq}

This subsection presents the notions of distribution and clustering of a matrix sequence, both in the sense of the eigenvalues and singular values, and some related definitions and results. In what follows, ${C}_c(\mathbb R)$ is the space of continuous complex-valued functions with bounded support on $\mathbb R$, and $ C_c(\mathbb C)$ is defined in the same way.

\begin{definition}\label{def-distribution}
Let $\{A_n\}_n$ be a matrix sequence, with $A_n$ of size
$d_n$, and let $\psi:D\subset\mathbb R^t\to\mathbb{C}^{r\times r}$ be a measurable function defined on a set $D$ with
$0<\mu_t(D)<\infty$.
\begin{itemize}
    \item We say that $\{A_n\}_n$ has a (asymptotic) singular value distribution described by $\psi$, and we write $\{A_n\}_n\sim_\sigma \psi$, if
    \begin{equation*} 
     \lim_{n\to\infty}\frac1{d_n}\sum_{i=1}^{d_n}F(\sigma_i(A_n))=\frac1{\mu_t(D)}\int_D\frac{\sum_{i=1}^{r}F(\sigma_i(\psi(\mathbf x)))}{r}{\rm d}\mathbf x,\qquad\forall\,F\in C_c(\mathbb R).
    \end{equation*}
    \item We say that $\{A_n\}_n$ has a (asymptotic) spectral (or eigenvalue) distribution described by $\psi$, and we write $\{A_n\}_n\sim_\lambda \psi$, if
    \begin{equation*}
     \lim_{n\to\infty}\frac1{d_n}\sum_{i=1}^{d_n}F(\lambda_i(A_n))=\frac1{\mu_t(D)}\int_D\frac{\sum_{i=1}^{r}F(\lambda_i(\psi(\mathbf x)))}{r}{\rm d}\mathbf x,\qquad\forall\,F\in C_c(\mathbb C).
    \end{equation*}
\end{itemize}
If $\psi$ describes both the singular value and eigenvalue distribution of $\{A_n\}_n$, we write $\{A_n\}_n\sim_{\sigma,\lambda}\psi$.
\end{definition}



The informal meaning behind the spectral distribution definition is the following: if $\psi$ is continuous, then a suitable ordering of the eigenvalues $\{\lambda_j(A_n)\}_{j=1,\ldots,d_n}$, assigned in correspondence with an equispaced grid on $D$, reconstructs approximately the $r$ surfaces $\mathbf x$ $\mapsto\lambda_i(\psi(\mathbf x)),\ i=1,\ldots,r$.
For instance, in the simplest case where $t=1$ and $D=[a,b]$, $d_n=nr$, the eigenvalues of $A_n$ are approximately equal - up to few potential outliers - to $\lambda_i\big(\psi(x_j)\big)$, where $x_j = a+j\frac{(b-a)}{n}$, $j=1,\ldots,n,\ i=1,\ldots,r$.
If $t=2$ and and $D=[a_1,b_1]\times [a_2,b_2]$, $d_n=n^2 r$, the eigenvalues of $A_n$ are approximately equal - again up to few potential outliers -  to $\lambda_i\big(\psi(x_{j_1},y_{j_2})\big)$, where $x_{j_1} = a_1+j_1\frac{b_1-a_1}{n}, y_{j_2}= a_2+j_2\frac{b_2-a_2}{n}$,  $j_1,j_2=1,\ldots,n,\ i=1,\ldots,r$. For $t\geq 3$ a similar reasoning applies.


Now let us proceed to the notion of clustering. We recall that the $\epsilon$-expansion of $S\subseteq\mathbb C$ is the set $B(S,\epsilon) := \bigcup_{z\in S}B(z,\epsilon)$, where $B(z,\epsilon) := \{w\in\mathbb C:\,|w-z|<\epsilon\}$ is the complex disk with center $z$ and radius $\epsilon>0$.

\begin{definition}
Let $\{A_n\}_n$ be a matrix sequence, with $A_n$ of size $d_n$, and let $S\subseteq\mathbb C$ be nonempty and closed. $\{A_n\}_n$ is {\em strongly clustered} at $S$ in the sense of the eigenvalues if, for any $\epsilon>0$ and as $n$ tends to infinity, the number of eigenvalues of $A_n$ outside $B(S,\epsilon)$ is bounded by a constant $q_\epsilon$ independent of $n$. In symbols,
$$q_\epsilon(n,S) := \#\big\{j\in\{1,\ldots,d_n\}: \lambda_j(A_n)\notin B(S,\epsilon)\big\}=O(1),\quad\mbox{as $n\to\infty$.}$$
$\{A_n\}_n$ is {\em weakly clustered} at $S$ if, for each $\epsilon>0$,
$$q_\epsilon(n,S)=o(d_n), \quad \mbox{as $n\to\infty.$}$$
A corresponding definition can be given for the singular values, with $S\subseteq\mathbb{R}^+$. If $\{A_n\}_n$ is strongly or weakly clustered at $S$ and $S$ is not connected, the connected components of $S$ are called sub-clusters.
\end{definition}

The case of spectral single point clustering, where $S$ consists of a single complex number $s$, retains special significance in the theory of preconditioning.


In the following remark, we clarify the deep connection between the two definitions above; in fact, clustering can be seen as a special case of distribution. First, recall that given a measurable function $g:D\subseteq\mathbb R^t\to\mathbb C$, the \emph{essential range} of $g$ is defined as $\mathcal{ER}(g) := \big\{z\in\mathbb C:\,\mu_t\big(\{g\in B(z,\epsilon)\}\big)>0\mbox{ for all $\epsilon>0$}\big\}$, where $\{g\in B(z,\epsilon)\} := \{x\in D:\,g(x)\in B(z,\epsilon)\}$. Generalizing the concept to a matrix-valued function $\psi:D\subset\mathbb R^t\to\mathbb{C}^{r\times r}$, the essential range is the union of the essential ranges of the eigenvalue functions $\lambda_i(\psi),\ i=1,\ldots, r$.

\begin{remark}
If $\{A_n\}_n\sim_\lambda \psi$, then $\{A_n\}_n$ is weakly clustered at the essential range of $\psi$ (see \cite[Theorem 4.2]{orth-distribution2}). Furthermore, if $\mathcal{ER}(\psi)=\{s\}$ with $s$ a fixed complex number, then $\{A_n\}_n\sim_\lambda \psi$ iff $\{A_n\}_n$ is weakly clustered at $s$ in the sense of the eigenvalues.

These statements can be translated to the singular value setting as well, with obvious minimal modifications. For instance, if $\mathcal{ER}(|\psi|)=s$ with $s$ a fixed nonnegative number, then $\{A_n\}_n\sim_\sigma \psi$ iff $\{A_n\}_n$ is weakly clustered at $s$ in the sense of the singular values.
\end{remark}

\subsection{Approximating classes of sequences}\label{ssez:acs}

In this subsection, we present the notion of approximating class of sequences and a related key result.

\begin{definition}{\rm (Approximating class of sequences)}\label{def:ACS}
	Let $\{A_n\}_{{n}}$ be a matrix-sequence and let $\{\{B_{n,j}\}_{{n}}\}_j$ be a class of matrix-sequences, with $A_n$ and $B_{n,j}$ of size $d_n$. We say that $\{\{B_{n,j}\}_n\}_j$ is an {\em approximating class of sequences (a.c.s.)} for $\{A_n\}_n$ if the following condition is met: for every $j$ there exists $n_j$ such that, for every $n\ge n_j$,
	\[ A_{{n}}=B_{{{n}},j}+R_{{{n}},j}+N_{{{n}},j}, \]
	\[ \textrm{rank}~R_{{{n}},j}\leq c(j) d_n  \quad {\rm and} \quad \|N_{{n},j}\|\leq\omega(j), \]
	where ${n}_j$, $c(j)$, and $\omega(j)$ depend only on $j$ and \[\lim_{j\to\infty}c(j)=\lim_{j\to\infty}\omega(j)=0.\]
 $\{\{B_{{n},j}\}_{{n}}\}_j \xrightarrow{\text{a.c.s.\ wrt\ $j$}}\{A_{{n}}\}_{{n}}$ denotes that $\{\{B_{{{n}},j}\}_{{n}}\}_j$ is an a.c.s. for $\{A_{{n}}\}_{{n}}$.
\end{definition}
	


The following theorem represents the expression of a related convergence theory and it is a powerful tool used, for example, in the construction of the GLT $*$-algebra.

\begin{theorem}\label{lem:Corollary5.1 and 5.2 in bookI}
Let $\{A_n\}_n$, $\{B_{n,j}\}_n$, with $j,n\in \mathbb{N}$, be matrix-sequences and let $\psi,\psi_j:D \subset \mathbb{R}^d \to \mathbb{C}$ be measurable functions defined on a set $D$ with positive and finite Lebesgue measure. Suppose that
\begin{enumerate}
	\item $\{B_{n,j}\}_n \sim_{\sigma} \psi_j$ for every $j$;
	\item $\{\{B_{n,j}\}_n\}_j \xrightarrow{\text{a.c.s.\ wrt\ $j$}} \{A_n\}_n$;
	\item $\psi_j \to \psi$ in measure.
\end{enumerate}
Then
\[ \{A_{{n}}\}_{{n}} \sim_{\sigma}  \psi. \]
Moreover, if all the involved matrices are Hermitian, the first assumption is replaced by $\{B_{{{n}},j}\}_{n}\sim_{\lambda}  \psi_j$ for every $j$, and the other two are left unchanged, then  $\{A_{{n}}\}_{{n}} \sim_{\lambda}  \psi$.
\end{theorem}


\subsection{Matrix structures}\label{ssez:Toeplitz-etc}

In this subsection we introduce the three types of matrix structures that constitute the basic building blocks of the GLT $*$-algebra.

\paragraph{Zero-distributed sequences.}
Zero-distributed sequences are defined as matrix sequences $\{A_n\}_n$ such that $\{A_n\}_n\sim_\sigma 0$. Note that, for any $r\ge1$, $\{A_n\}_n\sim_\sigma 0$ is equivalent to $\{A_n\}_n\sim_\sigma O_r$, where $O_r$ is the $r\times r$ zero matrix. The following theorem, taken from \cite{GLT-BookI,taud2}, provides a useful characterization for detecting this type of sequences. We use the natural convention $1/\infty=0$.

\begin{theorem}\label{0cs}
Let $\{A_n\}_n$ be a matrix sequence, with $A_n$ of size $d_n$. Then
\begin{itemize}
    \item $\{A_n\}_n\sim_\sigma 0$ if and only if $A_n=R_n+N_n$ with ${\rm rank}(R_n)/d_n\to0$ and $\|N_n\|\to0$ as $n\to\infty$;
    \item $\{A_n\}_n\sim_\sigma 0$ if there exists $p\in[1,\infty]$ such that $\|A_n\|_p/(d_n)^{1/p}\to0$ as $n\to\infty$.
\end{itemize}
\end{theorem}

\paragraph{Multilevel block Toeplitz matrices.}

Given $\mi{n}\in \mathbb{N}^d$, a matrix of the form
\[
[A_{\mi{i}-\mi{j}}]_{\mi{i},\mi{j}=\textbf{1}}^{\mi{n}} \in \mathbb{C}^{\nu(\mi{n})r \times \nu(\mi{n})r},
\]
with blocks $A_\mi{k} \in \mathbb{C}^{r\times r}$, $\mi{k} \in [-(\mi{n}-\textbf{1}), \ldots, \mi{n}-\textbf{1}]$, is called a multilevel block Toeplitz matrix, or, more precisely, a $d$-level $r$-block Toeplitz matrix.

Given a matrix-valued function $f : [-\pi, \pi]^{d} \rightarrow \mathbb{C}^{r\times r}$ belonging to $L^1([-\pi, \pi]^d)$, the $\mi{n}$-th Toeplitz matrix associated with $f$ is defined as
\begin{equation*} 
T_\mi{n}(f) := [\hat{f}_{\mi{i}-\mi{j}}]_{\mi{i},\mi{j}=\textbf{1}}^{\mi{n}} \in \mathbb{C}^{\nu(\mi{n})r \times \nu(\mi{n})r},
\end{equation*}
where
\begin{equation*} 
\hat{f}_\mi{k} =\frac{1}{(2\pi)^d} \int_{[-\pi,\pi]^d}
f(\boldsymbol{\theta})e^{-\iu (\mi{k}, \boldsymbol{\theta})}
{\rm d}\boldsymbol{\theta} \in \mathbb{C}^{r\times r}, \qquad \mi{k}
\in \mathbb{Z}^d,
\end{equation*}
are the Fourier coefficients of $f$, in which $\iu$ denotes the imaginary unit, the integrals are computed componentwise and $(\mi{k},\boldsymbol{\theta}) = k_1\theta_1 + \ldots + k_d\theta_d$. Equivalently, $T_\mi{n}(f)$ can be expressed as
\begin{equation*} 
T_\mi{n}(f) = \sum_{|j_1|<n_1} \ldots \sum_{|j_d|<n_d} \big[J_{n_1}^{(j_1)} \otimes \ldots \otimes J_{n_d}^{(j_d)}\big] \otimes \hat{f}_{(j_1,\ldots, j_d)}
\end{equation*}
where $\otimes$ denotes the Kronecker tensor product between matrices and $J_m^{(l)}$ is the matrix of order $m$ whose $(i,j)$ entry equals $1$ if $i-j=l$ and zero otherwise.

$\{T_\mi{n}(f)\}_{\mi{n}\in \mathbb{N}^d}$ is the family of (multilevel block) Toeplitz matrices associated with $f$, which is called the generating function.

\paragraph{Block diagonal sampling matrices.}


Given $d\ge 1$, $\mi{n}\in\mathbb{N}^d$ and a function $a:[0,1]^d\to\mathbb C^{r\times r}$, we define the multilevel block diagonal sampling matrix $D_\mi{n}(a)$ as the block diagonal matrix
\[
D_\mi{n}(a)=\diag_{\mi{i}=\mathbf{1},\ldots,\mi{n}}a\Bigl(\frac{\mi{i}}{\mi{n}}\Bigr)\in\mathbb C^{\nu(\mi{n})r \times \nu(\mi{n})r}.
\]

\subsection{The $*$-algebra of multilevel block GLT sequences}\label{ssez:GLT}
Let $r\ge1$ be a fixed integer. A multilevel $r$-block GLT sequence, or simply a GLT sequence if we do not need to specify $r$, is a special multilevel $r$-block matrix-sequence equipped with a measurable function $\kappa:[0,1]^d\times[-\pi,\pi]^d\to\mathbb C^{r\times r}$, $d\ge 1$, called \emph{symbol}. The symbol is essentially unique, in the sense that if $\kappa$, $\varsigma$ are two symbols of the same GLT sequence, then $\kappa=\varsigma$ a.e.. We write $\{A_n\}_n\sim_{\rm GLT}\kappa$ to denote that $\{A_n\}_n$ is a GLT sequence with symbol $\kappa$.


It can be proven that the set of multilevel block GLT sequences is the $*$-algebra generated by the three classes of sequences defined in Subsection \ref{ssez:Toeplitz-etc}: zero-distributed, multilevel block Toeplitz and block diagonal sampling matrix sequences. The GLT class satisfies several algebraic and topological properties that are treated in detail in \cite{GLT-block1D,GLT-blockdD,GLT-BookI,GLT-BookII}. Here, we focus on the main operative properties, listed below, that represent a complete characterization of GLT sequences, equivalent to the full constructive definition.

\begin{enumerate}
    \item[\textbf{GLT\,1.}] If $\{A_n\}_n\sim_{\rm GLT}\kappa$ then $\{A_n\}_n\sim_\sigma\kappa$ in the sense of Definition \ref{def-distribution}, with $t=2d$ and $D=[0,1]^d\times[-\pi,\pi]^d$. If moreover each $A_n$ is Hermitian, then $\{A_n\}_n\sim_\lambda\kappa$, again in the sense of Definition \ref{def-distribution} with $t=2d$.
    \item[\textbf{GLT\,2.}] We have:
    \begin{itemize}
        \item $\{T_\mi{n}(f)\}_\mi{n}\sim_{\rm GLT}\kappa(\mathbf{x},\boldsymbol{\theta})=f(\boldsymbol{\theta})$ if $f:[-\pi,\pi]^d\to\mathbb C^{r\times r}$ is in $L^1([-\pi,\pi]^d)$;
        \item $\{D_\mi{n}(a)\}_\mi{n}\sim_{\rm GLT}\kappa(\mathbf{x},\boldsymbol{\theta})=a(\mathbf{x})$ if $a:[0,1]^d\to\mathbb C^{r\times r}$ is Riemann-integrable;
        \item $\{Z_n\}_n\sim_{\rm GLT}\kappa(\mathbf{x},\boldsymbol{\theta})=O_r$ if and only if $\{Z_n\}_n\sim_\sigma 0$.
    \end{itemize}
    \item[\textbf{GLT\,3.}] If $\{A_n\}_n\sim_{\rm GLT}\kappa$ and $\{B_n\}_n\sim_{\rm GLT}\varsigma$ then:
    \begin{itemize}
        \item $\{A_n^*\}_n\sim_{\rm GLT}\kappa^*$;
        \item $\{\alpha A_n+\beta B_n\}_n\sim_{\rm GLT}\alpha\kappa+\beta\varsigma$ for all $\alpha,\beta\in\mathbb C$;
        \item $\{A_nB_n\}_n\sim_{\rm GLT}\kappa\varsigma$;
        \item $\{A_n^\dag\}_n\sim_{\rm GLT}\kappa^{-1}$ provided that $\kappa$ is invertible a.e.
    \end{itemize}
    \item[\textbf{GLT\,4.}] $\{A_n\}_n\sim_{\rm GLT}\kappa$ if and only if there exist $\{B_{n,j}\}_n\sim_{\rm GLT}\kappa_j$ such that $\{\{B_{n,j}\}_n\}_j \xrightarrow{\text{a.c.s.\ wrt\ $j$}} \{A_n\}_n$ and $\kappa_j\to\kappa$ in measure.
\end{enumerate}

Note that, by GLT 1, it is always possible to obtain the singular value distribution from the GLT symbol, while the eigenvalue distribution can be deduced only if the involved matrices are Hermitian. However, interesting tools are available for matrix sequences in which the non-Hermitian part is somehow negligible, see \cite{orth-distribution2,orth-distribution1} and references therein. The main result is reported below, more advanced tools can be found in \cite{non-herm-perturb}.

\begin{theorem}{\cite[Corollary 4.1]{GLT-BookI}\cite[Corollary 2.3]{GLT-BookII}} \label{th:GLT-lambda}
Let $\{X_n\}_n,\,\{Z_n\}_n$ be matrix sequences, with $X_n$, $Z_n$ of size $d_n$, and set $Y_n=X_n+Z_n$. Assume that the following conditions are met
\begin{enumerate}
	\item $\|X_n\|,\|Z_n\|\le C$ for all $n$, where $C$ is a constant independent of $n$;
  \item Every $X_n$ is Hermitian and $\{X_n\}_n\sim_\lambda \kappa$;
  \item $\|Z_n\|_1=o(d_n)$.
\end{enumerate}
Then $\{Y_n\}_n\sim_\lambda \kappa$. Moreover, if $\|Z_n\|_1=O(1)$, the range of $\kappa$ is a strong cluster for the eigenvalues of $\{Y_n\}_n$.
\end{theorem}

\section{Discretized problem and GLT analysis}\label{sec:appr}

In this section, we briefly present the discretization of the two-dimensional fractional Helmholtz equation \eqref{a1}, with $\Omega=[0,1]^{2}$ and $\Omega^{c} = \mathbb{R}^{2}\setminus\Omega$, and then compute the GLT symbol of the resulting coefficient matrix sequence in order to study the singular value and eigenvalue distribution.

\subsection{Discretization}

We follow the fractional centered differences techique adopted in \cite{20}, to which we refer for further details.

Let $n$ be a positive integer and define the uniform spatial partition $\Omega_{h}=\{(x_{i},y_{j}) : i,j=1,2,\ldots n\}$, with stepsize $h=\frac{1}{n+1}$ and grid points $x_{i}=ih$, $y_{j}=jh$, for $i,j\in \mathbb{Z}$. Consider the discrete fractional Laplacian operator
\begin{equation}\label{b1}
    (-\Delta_{h})^{\alpha/2} u(x,y) :=
    \frac{1}{h^{\alpha}}\sum_{k_{1},k_{2}\in\mathbb{Z}}b_{k_{1},k_{2}}^{(\alpha)}u(x+k_{1}h,y+k_{2}h),
\end{equation}
where
\begin{equation*}
    b_{k_{1},k_{2}}^{(\alpha)}=\frac{1}{4\pi^{2}}\int_{-\pi}^{\pi}\int_{-\pi}^{\pi}t_\alpha (\theta_1,\theta_2) e^{-\iu(k_{1}\theta_1+k_{2}\theta_2)}{\rm d}\theta_1 {\rm d}\theta_2, \qquad \iu^2 = -1,
\end{equation*}
are the Fourier coefficients of the function
\begin{equation*} 
    t_\alpha (\theta_1,\theta_2) = \left[4\sin^2\left(\frac{\theta_1}{2}\right)+4\sin^2\left(\frac{\theta_2}{2}\right)\right]^{\frac{\alpha}{2}}.
\end{equation*}
Now define the weighted space
\begin{equation*}
    W^{\gamma,1}(\mathbb{R}^{2})=\left\{ u\in L^{1}(\mathbb{R}^{2}) : \int_{\mathbb{R}^{2}}(1+|z|)^{\gamma}|\Tilde{u}(\Tilde{x},\Tilde{y})| \,{\rm d}\Tilde{x}{\rm d}\Tilde{y} < \infty \right \},
\end{equation*}
in which $|z|^{2}=\Tilde{x}^{2}+\Tilde{y}^{2}$ and $\Tilde{u}(\Tilde{x},\Tilde{y})$ is the Fourier transform of $u(x,y)$. Thanks to \cite[Theorem 2.2]{20}, for any $u\in W^{2+\alpha,1}(\mathbb{R}^{2})$ it holds
\begin{equation}\label{trunc-err}
    (-\Delta)^{\alpha/2}u(x,y)=(-\Delta_{h})^{\alpha/2}u(x,y) + \mathcal{O}(h^2).
\end{equation}
Combining \eqref{b1} and \eqref{trunc-err} and applying them to \eqref{a1}, we obtain
\begin{equation*}
    (-\Delta_{h})^{\alpha/2}u(x_{i},y_{j})+\mu(x_i,y_j)u(x_{i},y_{j})=v(x_{i},y_{j})+\Tilde{R}_{i,j,h}, \qquad (x_{i},y_{j})\in \Omega_{h},
\end{equation*}
where $\Tilde{R}_{i,j,h}$ is the truncation error. Since the approximation is uniform, there exists a constant $c_u$ independent of $i,j,h$ such that $|\Tilde{R}_{i,j,h}|\leq c_{u}h^{2}$, therefore the method is consistent with precision order $2$.

Letting $\mu_{i,j}=\mu(x_{i},y_{j})$, $v_{i,j}=v(x_{i},y_{j})$, $u_{i,j}\approx u(x_i,y_j)$ and disregarding the truncation error, we arrive at the following numerical scheme for solving the continuous problem
  \begin{equation*} 
      \begin{cases}
           (-\Delta_h)^{\alpha/2} u_{i,j}+\mu_{i,j}u_{i,j}=v_{i,j}, & (x_i,y_j)\in\Omega_h, \\
           u_{i,j}=0, & (x_i,y_j)\in\Omega_h^c,
      \end{cases}
  \end{equation*}
where $\Omega_{h}^{c} := \{(x_{i},y_{j}) : i,j\in\mathbb{Z}\}\setminus\Omega_{h}$.

Making use of the multi-index notation with the standard lexicographic order and setting
\begin{align*}
    \vet{u} &= \begin{bmatrix} u_{1,1} & u_{1,2} & \ldots & u_{1,n} & u_{2,1} & u_{2,2} & \ldots & u_{2,n} & \ldots & u_{n,n} \end{bmatrix}, \\
    \vet{v} &= \begin{bmatrix} v_{1,1} & v_{1,2} & \ldots & v_{1,n} & v_{2,1} & v_{2,2} & \ldots & v_{2,n} & \ldots & v_{n,n} \end{bmatrix},
\end{align*}
the scheme above is rewritten as the linear system
\begin{equation*}
    A_{\mi{n}} \vet{u}:=\big(B_{\mi{n}}+D_{\mi{n}}(\mu)\big) \vet{u} = \vet{v}, \qquad {\mi{n}}=(n,n)
\end{equation*}
where $B_{\mi{n}}=\frac{1}{h^{\alpha}}T_{\mi{n}}(t_\alpha)$ and $D_{\mi{n}}(\mu)$ is the block diagonal sampling matrix associated to $\mu$. More explicitly, the two-level symmetric Toeplitz matrix $T_{\mi{n}}(t_\alpha)$ generated by $t_\alpha(\theta_1,\theta_2)$ has the form
\begin{equation*}
T_{\mi{n}}(t_\alpha)=\begin{bmatrix}
    C_{0} & C_{1} & C_{2} & \cdots & C_{n-2} & C_{n-1}\\
    C_{1} & C_{0} & C_{1} & \cdots & C_{n-3} & C_{n-2} \\
    C_{2} & C_{1} & C_{0} & \cdots & C_{n-4} & C_{n-3} \\
    \vdots & \vdots &\vdots& \ddots & \vdots & \vdots\\
    C_{n-2} & C_{n-3} & C_{n-4} & \cdots & C_{0} & C_{1}\\
    C_{n-1} & C_{n-2} & C_{n-3} & \cdots & C_{1} & C_{0}
\end{bmatrix},
\end{equation*}
in which
\begin{equation*}
    C_{k} = \begin{bmatrix}
        b_{0,k}^{(\alpha)} & b_{1,k}^{(\alpha)} & b_{2,k}^{(\alpha)} & \cdots & b_{n-2,k}^{(\alpha)} & b_{n-1,k}^{(\alpha)}\\
        b_{1,k}^{(\alpha)} & b_{0,k}^{(\alpha)} & b_{1,k}^{(\alpha)} & \cdots & b_{n-3,k}^{(\alpha)} & b_{n-2,k}^{(\alpha)} \\
        b_{2,k}^{(\alpha)} & b_{1,k_2}^{(\alpha)} & b_{0,k}^{(\alpha)} & \cdots & b_{n-4,k}^{(\alpha)} & b_{n-3,k}^{(\alpha)} \\
        \vdots & \vdots &\vdots& \ddots & \vdots & \vdots\\
        b_{n-2,k}^{(\alpha)} & b_{n-3,k}^{(\alpha)} & b_{n-4,k}^{(\alpha)} & \cdots & b_{0,k}^{(\alpha)} & b_{1,k}^{(\alpha)}\\
        b_{n-1,k}^{(\alpha)} & b_{n-2,k}^{(\alpha)} & b_{n-3,k}^{(\alpha)} & \cdots & b_{1,k}^{(\alpha)} & b_{0,k}^{(\alpha)}
    \end{bmatrix}.
\end{equation*}
\begin{remark}
    The coefficients $b_{k_1,k_2}^{(\alpha)}$ can be computed via the trapezoidal rule
    \begin{align*}
        b_{k_1,k_2}^{(\alpha)}
        &\approx \frac{1}{M^{2}}\sum_{p=0}^{M-1}\sum_{q=0}^{M-1}t_\alpha(p\delta,q\delta) e^{-\iu(k_1 p \delta+ k_2 q\delta)},
        \qquad \delta=\frac{2\pi}{M}, \quad M\gg 2n.
    \end{align*}
    Using the FFT, the computation of $b_{k_1,k_2}^{(\alpha)}$, $ k_1,k_2 = 0,\ldots, M-1$, costs $\mathcal{O}(M \log M)$, then in practice only the coefficients with $0\leq k_1,k_2 \leq n-1$ are employed.
\end{remark}

\subsection{GLT, singular value and eigenvalue analysis}

This subsection contains the GLT analysis of the coefficient matrix sequence obtained previously. In order to obtain a meaningful distribution, we rescale the linear system as
    \begin{equation}\label{scaled eq}
     \hat{A}_{\mi{n}} \vet{u} := \big(\hat{B}_{\mi{n}}+h^\alpha D_{\mi{n}}(\mu)\big) \vet{u} = h^\alpha \vet{v}, \qquad \mi{n}=(n,n),
  \end{equation}
in which $\hat{B}_{\mi{n}} := T_{\mi{n}}(t_\alpha)$.

\begin{remark} \label{remark:toplitz-part}
It holds
\begin{enumerate}
    \item $\{\hat{B}_{\mi{n}}\}_n\sim_{\rm GLT} t_\alpha$;
    \item $\{\hat{B}_{\mi{n}}\}_n\sim_\sigma t_\alpha$;
    \item $\{\hat{B}_{\mi{n}}\}_n\sim_\lambda t_\alpha$.
\end{enumerate}
In fact, the first assertion is just an application of the first item of {\bf GLT 2}. The second and third one are consequences of {\bf GLT 1}, since $\hat{B}_{\mi{n}}$ is Hermitian.
\end{remark}

First, we perform the analysis with the assumption that the wave function $\mu(x,y)$ is essentially bounded. In this setting, we are able to obtain both the precise eigenvalue and spectral distribution. The case of an unbounded $\mu(x,y)$ is more delicate and is treated afterwards.

\subsubsection{The case of an essentially bounded wave function} \label{ssez:wave-function-bounded}

To begin, we observe that from the fact that $\mu(x,y)$ is essentially bounded it follows the inequality
\begin{equation} \label{eq:inequality-diagonal-part-bounded}
    \|h^\alpha D_{\mi{n}}(\mu)\| \le h^\alpha \|\mu\|_\infty.
\end{equation}
As a consequence, the following lemma can be proved. It will allow us to prove the first main distributional theorem.
\begin{lemma} \label{lemma:diagonal-part-bounded}
Let $\mu(x,y)$ be an essentially bounded complex-valued function. Then
\begin{enumerate}
    \item $\{h^\alpha D_{\mi{n}}(\mu)\}_n\sim_{\rm GLT} 0$;
    \item $\{h^\alpha D_{\mi{n}}(\mu)\}_n\sim_\sigma 0$;
    \item $\{h^\alpha D_{\mi{n}}(\mu)\}_n\sim_\lambda 0$.
\end{enumerate}
\end{lemma}
\begin{prova}
Points 2 and 3 follow directly from \eqref{eq:inequality-diagonal-part-bounded}:
\begin{enumerate}
    \item[2.] all the singular values of $h^\alpha D_{\mi{n}}(\mu)$ converge to zero, at least as fast as $h^\alpha \|\mu\|_\infty$. Then the related sequence $\{h^\alpha D_{\mathbf{n}}(\mu)\}_n$ is strongly clustered at 0 in the sense of the singular values and it follows from Theorem \ref{0cs} that $\{h^\alpha D_{\mathbf{n}}(\mu)\}_n$ is zero-distributed; either by the first item, taking $R_n$ equal to the null matrix, or the second one, taking $p=\infty$;
    \item[3.] all the eigenvalues of $h^\alpha D_{\mi{n}}(\mu)$ converge to zero, at least as fast as $h^\alpha \|\mu\|_\infty$, since the spectral radius of a matrix is bounded from above by its spectral norm. By direct inspection, $\{h^\alpha D_{\mi{n}}(\mu)\}_n\sim_\lambda 0$.
\end{enumerate}
Finally, 1 follows from 2 simply by the third item in {\bf GLT 2}.
\end{prova}

\begin{theorem}\label{th:coeff matrix-sequence SV}
Let $\mu(x,y)$ be an essentially bounded complex-valued function.
Then
\begin{itemize}
    \item[\bf a1.] $\{\hat{B}_{\mi{n}}+h^\alpha D_{\mi{n}}(\mu)\}_n \sim_{\rm GLT} t_\alpha$;
    \item[\bf a2.] $\{\hat{B}_{\mi{n}}+h^\alpha D_{\mi{n}}(\mu)\}_n\sim_\sigma t_\alpha$;
    \item[\bf a3.] $\{\hat{B}_{\mi{n}}+h^\alpha D_{\mi{n}}(\mu)\}_n\sim_\lambda t_\alpha$.
\end{itemize}
\end{theorem}

\begin{prova}
Point {\bf a1} follows from the first items of Remark \ref{remark:toplitz-part} and Lemma \ref{lemma:diagonal-part-bounded} with a straightforward application of the GLT tools: invoking the $*$-algebra structure of the GLT class expressed by {\bf GLT 3}, we obtain that $\{\hat{B}_{\mi{n}}+h^\alpha D_{\mi{n}}(\mu)\}_n\sim_{\rm GLT} t_\alpha + 0 =t_\alpha$. Then {\bf a2} is an obvious consequence of {\bf a1} and {\bf GLT 1}.

Regarding {\bf a3}, we invoke Theorem \ref{th:GLT-lambda} with $X_n = \hat{B}_{\mi{n}}$ and $Z_n = h^\alpha D_{\mi{n}}(\mu)$. Since $\hat{B}_{\mi{n}}$ is Hermitian, we can apply known results concerning the extreme eigenvalues of Hermitian Toeplitz matrices \cite{S-extr2,S-extr1} and obtain that $\|\hat{B}_{\mi{n}}\|=\|T_{\mi{n}}(t_\alpha)\| < \|t_\alpha\|_\infty$, where the strict inequality holds because $\min t_\alpha = 0$ and $\max t_\alpha = 2^{3\alpha/2} > 0$, in accordance with Theorem 2.1 in \cite{S-extr1} by setting $g=1$ and $f=t_\alpha$. Recalling \eqref{eq:inequality-diagonal-part-bounded}, the first assumption of Theorem \ref{th:GLT-lambda} is satisfied. Moreover, the second assumption holds true by Remark \ref{remark:toplitz-part}. Finally, inequality \eqref{eq:inequality-diagonal-part-bounded} leads to
\[
\|h^\alpha D_{\mi{n}}(\mu)\|_1 \le h^\alpha \|\mu\|_\infty n^2=o(n^2),
\]
which is the third and last assumption. Therefore by Theorem \ref{th:GLT-lambda} we can conclude $\{\hat{B}_{\mi{n}}+h^\alpha D_{\mi{n}}(\mu)\}_n\sim_\lambda t_\alpha$.
\end{prova}

\subsubsection{The case of an unbounded wave function}

Now we proceed to the case where $\mu(x,y)$ is unbounded. For this purpose, for any $M>0$ we define the following truncated version of $\mu$
\begin{equation}\label{mu truncated}
\mu_M(x,y)=
\left\{
\begin{array}{cl}
\mu(x,y) & {\rm if}\, \, |\mu(x,y)|\le M, \\
M        & {\rm otherwise,}
\end{array}
\right.
\end{equation}
to which, in analogy with the structure of $D_{\mi{n}}(\mu)$, we associate the block diagonal sampling matrix $D_{\mi{n}}(\mu_M)$. This process constructs a class of matrix sequences $\{\{h^\alpha D_{\mi{n}}(\mu_M)\}_n\}_M$, allowing us to rely on the notion of a.c.s. in the sense of Definition \ref{def:ACS}. We begin with a preliminary lemma.

\begin{lemma} \label{lemma:diagonal-part-unbounded}
Let $\mu_M(x,y)$ be Riemann-integrable for any $M>0$ and suppose that $\{\mu_M(x,y)\}_M$ converges pointwise to $\mu(x,y)$, with $\mu(x,y)$ unbounded and complex-valued. Then
\begin{enumerate}
    \item $\{h^\alpha D_{\mi{n}}(\mu)\}_n\sim_{\rm GLT} 0$;
    \item $\{h^\alpha D_{\mi{n}}(\mu)\}_n\sim_\sigma 0$.
\end{enumerate}
\end{lemma}
\begin{prova}
    Let us start with the proof of 2. From the hypothesis that $\mu_M$ is Riemann-integrable for any $M>0$ and the second item of {\bf GLT 2}, we deduce that $\{D_{\mi{n}}(\mu_M)\}_n\sim_{\rm GLT}\mu_M$. Now consider the sequence $\{h^\alpha I_{n^2}\}_n$: it is obviously zero-distributed, since any singular value equals $h^\alpha$, therefore it is a GLT sequence with symbol 0 by the third item of {\bf GLT 2}. So, by the third item of {\bf GLT 3}, we obtain
\[
\{h^\alpha D_{\mi{n}}(\mu_M) = h^\alpha I_{n^2} \cdot D_{\mi{n}}(\mu_M) \}_n \sim_{\rm GLT} 0 \cdot \mu_M = 0,
\]
or equivalently $\{h^\alpha D_{\mi{n}}(\mu_M)\}_n$ is zero-distributed. It is easy to see by direct inspection that $\{\{h^\alpha D_{\mi{n}}(\mu_M)\}_n\}_M$ is an a.c.s for $\{h^\alpha D_{\mi{n}}(\mu)\}_n$, so we conclude by applying Theorem \ref{lem:Corollary5.1 and 5.2 in bookI} with $j=M$, $\psi_M=0$, and $\psi=0$, and by observing that the pointwise convergence implies the weaker type of convergence in measure.
Then the first assertion is just an obvious consequence of 2 by {\bf GLT 2}.
\end{prova}

\begin{theorem}\label{th:coeff matrix-sequence SV-EXT}
Let $\mu_M(x,y)$ be Riemann-integrable for any $M>0$ and suppose that $\{\mu_M(x,y)\}_M$ converges pointwise to $\mu(x,y)$, with $\mu(x,y)$ unbounded and complex-valued. Then
\begin{itemize}
    \item[\bf b1.] $\{\hat{B}_{\mi{n}}+h^\alpha D_{\mi{n}}(\mu)\}_n\sim_{\rm GLT} t_\alpha$;
    \item[\bf b2.] $\{\hat{B}_{\mi{n}}+h^\alpha D_{\mi{n}}(\mu)\}_n\sim_\sigma t_\alpha$.
\end{itemize}
\end{theorem}

\begin{prova}
The proof of {\bf b1} is once again a straightforward application of the GLT axioms: we know from Remark \ref{remark:toplitz-part} that $\{\hat{B}_{\mathbf{n}}\}_n \sim_{\rm GLT} t_\alpha$, combining this result with Lemma \ref{lemma:diagonal-part-unbounded} and the second item of {\bf GLT 3} we get $\{\hat{B}_{\mi{n}}+h^\alpha D_{\mi{n}}(\mu)\}_n\sim_{\rm GLT} t_\alpha$. As usual, with {\bf GLT 1} we obtain the singular value distribution in {\bf b2}.
\end{prova}

Regarding the eigenvalue distribution in the case of an unbounded wave function $\mu(x,y)$, so far we have been unable to obtain a satisfactory generalization for the eigenvalues as the one we have for the singular values, even using the more sophisticated tools in \cite{non-herm-perturb}. Theorem \ref{th:coeff matrix-sequence SV-EXT} represents only a partial answer. Hence, this question needs to be explored further in the future.

\begin{remark}
    Note that, if $\mu(x,y)$ is real-valued, then all the distributional results found in Subsection \ref{ssez:wave-function-bounded} hold also in the unbounded case. In fact, under this assumption the diagonal part $h^\alpha D_{\mi{n}}(\mu)$ is Hermitian, so the eigenvalue distribution follows simply as a consequence of {\bf GLT 1}.
\end{remark}

\section{Preconditioning}\label{sec:sp-analysis}

In this section, after a brief introduction to the unilevel and two-level $\tau$ algebra, we consider the ad hoc $\tau$ preconditioners proposed in \cite{tau-prec} and perform the distribution and clustering analysis of the preconditioned matrix-sequences, both in the sense of the eigenvalues and singular values.

\subsection{The $\tau$ algebra} 

Given a symmetric Toeplitz matrix $\textit{T}_{n}\in \mathbb{R}^{{n}\times n}$, the natural $\tau$ preconditioner of $T_n$ is defined as
\begin{equation}\label{Hankel correction}
    \tau(T_n):= T_n-H(T_n),
\end{equation}
in which $H(T_{n})$ denotes the Hankel matrix, namely a matrix whose entries are constant along each antidiagonal, fully characterized by its first row $[t_{2},t_{3},\ldots,t_{n-1},0,0]$ and its last column $[0,0,t_{n-1}, \ldots,t_{3},t_{2}]^{\top}$, where $[t_{1},t_{2},\ldots,t_{n}]^{\top}$ is the first column of $T_n$.

A great amount of theoretical and computational research has been dedicated in the last decades to this type of preconditioner \cite{BC,CN,tau-theory1}. It is well-known that $\tau$ matrices form an algebra, closed under inversion, and that they can be uniformly diagonalized through the following discrete sine transform (DST) matrix
\begin{equation*}
    [S_n]_{i,j}=\sqrt{\frac{2}{n+1}}\sin{\left(\frac{\pi ij}{n+1}\right)}, \qquad 1\leq i, j\leq n.
\end{equation*}
In other words, we can always decompose $\tau(T_n)$ as
\begin{equation*} 
    \tau(T_{n})=S_n\Lambda_{n} S_n,
\end{equation*}
where $\Lambda_{n}$ is the diagonal matrix containing the eigenvalues of $\tau(\textit{T}_{n})$. Note that $S_n$ is real, symmetric and orthogonal, so that $S_n=S_n^T=S_n^{-1}$. Moreover, $S_n$ is associated to the fast sine transform of type I, thanks to which the DST (the product between $S_n$ and any vector in $\mathbb{R}^n$) can be performed in $\mathcal{O}(n\log{}n)$ real operations; more precisely the cost is about half of the renowned fast Fourier transform \cite{Van} (further details and a description of several sine/cosine transforms can be found in \cite{DiB-S,KO}). As a consequence, in the $\tau$ algebra all the relevant matrix operations, such as matrix-vector and matrix-matrix product, inversion, solution of a linear system or computation of the spectrum, have the computational cost of $\mathcal{O}(n\log{}n)$. In particular, the eigenvalues of $\tau(T_n)$ are obtained as the DST of the first column of $T_n$.

Generalizing to any $d$-level symmetric Toeplitz matrix $T_\mi{n}$, with $d\geq 1$ and $\mi{n}=(n_1,\ldots,n_d)$, the $d$-level $\tau$ preconditioner is diagonalized as
\begin{equation*}
    \tau(T_{\mi{n}})=S_\mi{n}\Lambda_{\mi{n}} S_\mi{n}, \qquad S_\mi{n} = S_{n_1} \otimes \cdots \otimes S_{n_d},
\end{equation*}
in which $\Lambda_{\mi{n}}$ is the diagonal matrix containing the eigenvalues of $\tau (T_\mi{n})$ and is obtained as the $d$-level discrete sine transform of type I of the first column of $T_{\mi{n}}$.

The set of $d$-level $\tau$ matrices is the $d$-level $\tau$ algebra. Again, all the relevant matrix operations in the $d$-level $\tau$ algebra cost $\mathcal{O}\big(\nu(\mi{n})\log \nu(\mi{n})\big)$ real operations, which, given the fact that the size of the matrices is $\nu(\mi{n})$, is quasi optimal.

From the algebraic point of view, the explicit construction of $\tau(T_\mi{n})$ can be done recursively: starting from the most external level and working inward, apply the additive decomposition (\ref{Hankel correction}) until the scalar level, where the matrices are $1 \times 1$, is reached. From distribution results on multilevel Hankel matrix sequences \cite{FaTi}, it is known for any $L^1$ function $f$ it holds
\begin{equation}\label{hankel:distr}
\left\{T_{\mi{n}}(f) - \tau\big(T_{\mi{n}}(f)\big)\right\}_n\sim_{\sigma,\lambda} 0.
\end{equation}

\subsection{Preconditioning proposals}

Here we give the definitions of two pairs of $\tau$ preconditioners that were proposed in \cite{tau-prec} for linear system \eqref{scaled eq}, in light of the excellent structural, spectral and computational features of the $\tau$ algebra.

The first pair is based on approximating $t_\alpha$, the function that generates the Toeplitz sequence in \eqref{scaled eq}, in a separable way as
\begin{equation*}
    \text{g}_{\alpha}^{(2)}(\theta_1,\theta_2)  = \text{g}_{\alpha}^{(1)}(\theta_1) + \text{g}_{\alpha}^{(1)}(\theta_2) = \bigg[4\sin^{2}\left(\frac{\theta_1}{2}\right)\bigg]^{\frac{\alpha}{2}}+\bigg[4\sin^{2}\left(\frac{\theta_2}{2}\right)
    \bigg]^{\frac{\alpha}{2}},
\end{equation*}
in which $\text{g}_{\alpha}^{(1)} (\theta) = \left[ 4\sin^{2} \left( \frac{\theta}{2} \right) \right]^{\frac{\alpha}{2}}$. Then, the 2-level Toeplitz matrix generated by $\text{g}_{\alpha}^{(2)}$ has the form
\begin{equation*}
    T_{\mi{n}}(\text{g}_{\alpha}^{(2)}) = I_n\otimes T_n(\text{g}_{\alpha}^{(1)}) + T_n(\text{g}_{\alpha}^{(1)})\otimes I_n, \qquad \mi{n}=(n,n).
\end{equation*}
Note that, combining {\bf GLT 2} and {\bf GLT 1} with the fact that $T_{\mi{n}}(\text{\rm g}_{\alpha}^{(2)})$ is real and symmetric (hence Hermitian), we have
\begin{equation}\label{remark:toeplitz-part-approx}
\{T_{\mi{n}}(\text{\rm g}_{\alpha}^{(2)})\}_n\sim_{\text{\rm GLT},\sigma,\lambda} \text{\rm g}_{\alpha}^{(2)}.
\end{equation}
The preconditioners are then defined as
\begin{equation*}
    \tau\big(T_{\mi{n}}(\text{g}_{\alpha}^{(2)})\big), \qquad
    \tau\big(T_{\mi{n}}(\text{g}_{\alpha}^{(2)})\big)+ h^\alpha \mu_h I_{n^2},
\end{equation*}
where $\mu_h$ is the arithmetic average of the evaluations of $\mu(x,y)$ over the grid points. In the case of a bounded variable coefficient $\mu(x,y)$, the authors in \cite{tau-prec} proved good localization results for the preconditioned matrices, with bounds independent of $n$. However, a straightforward GLT analysis demonstrates that the clustering at $1$ cannot hold.

\begin{prop} \label{prop:negative-preconditioner-1}
   It holds
    \begin{equation*}
        \left\{\tau\big(T_{\mi{n}}(\text{\rm g}_{\alpha}^{(2)})\big)^{-1}\left(\hat{B}_{\mi{n}}+h^\alpha D_{\mi{n}}(\mu)\right)\right\}_n\sim_{{\rm GLT},\sigma} \frac{t_\alpha}{\text{\rm g}_{\alpha}^{(2)}}.
    \end{equation*}
    Moreover, assuming that $\mu(x,y)$ is essentially bounded, it holds
    \begin{equation*}
        \left\{\tau\big(T_{\mi{n}}(\text{\rm g}_{\alpha}^{(2)})\big)^{-1}\left(\hat{B}_{\mi{n}}+h^\alpha D_{\mi{n}}(\mu)\right)\right\}_n\sim_{\lambda} \frac{t_\alpha}{\text{\rm g}_{\alpha}^{(2)}}.
    \end{equation*}
\end{prop}

\begin{prova}
     From Equation \eqref{hankel:distr} we know that $\big\{T_{\mi{n}}(\text{g}_{\alpha}^{(2)}) - \tau\big(T_{\mi{n}}(\text{g}_{\alpha}^{(2)})\big)\big\}_n$ is a zero-distributed sequence, hence it is a GLT sequence with symbol 0 by {\bf GLT 2}. Turning to the $*$-algebra properties in {\bf GLT 3} and recalling \eqref{remark:toeplitz-part-approx}, we deduce that $\big\{\tau(T_{\mi{n}}\big(\text{g}_{\alpha}^{(2)})\big)\big\}_n$ is a GLT sequence with symbol $\text{g}_{\alpha}^{(2)}$. Recalling that $\tau$ matrices are real and symmetric, using {\bf GLT 1}, {\bf GLT 3} and Theorems \ref{th:coeff matrix-sequence SV} and \ref{th:coeff matrix-sequence SV-EXT}, we conclude.
\end{prova}

\noindent Therefore the preconditioned sequence has eigenvalues clustered at the range of
${t_\alpha}/{\text{g}_{\alpha}^{(2)}}$, which is a nontrivial real positive interval and cannot reduce to the point $s=1$.

Before proceeding to the next preconditioner, we note that, very similarly to Lemma \ref{lemma:diagonal-part-bounded}, in the case where $\mu(x,y)$ is essentially bounded the maximal singular value of $h^\alpha \mu_h I_{n^2}$ is bounded from above by $h^\alpha \|\mu\|_\infty$, hence $\{h^\alpha \mu_h I_{n^2}\}_n$ is a zero-distributed GLT sequence. In other words, we have that
\begin{equation}\label{eq:diagonal-part-precond-distr}
    \{h^\alpha \mu_h I_{n^2}\}_n \sim_{{\rm GLT},\sigma} 0.
\end{equation}

\begin{prop}
    Assuming that $\mu(x,y)$ is essentially bounded, it holds
    \begin{equation*}
        \left\{\left(\tau\big(T_{\mi{n}}(\text{\rm g}_{\alpha}^{(2)})\big)+ h^\alpha {\mu}_h I_{n^2}\right)^{-1}\left(\hat{B}_{\mi{n}}+h^\alpha D_{\mi{n}}(\mu)\right)\right\}_n \sim_{{\rm GLT},\sigma} \frac{t_\alpha}{\text{\rm g}_{\alpha}^{(2)}}.
    \end{equation*}
    Moreover, if $\mu(x,y)$ is real-valued, it holds
    \begin{equation*}
        \left\{\left(\tau\big(T_{\mi{n}}(\text{\rm g}_{\alpha}^{(2)})\big)+ h^\alpha {\mu}_h I_{n^2}\right)^{-1}\left(\hat{B}_{\mi{n}}+h^\alpha D_{\mi{n}}(\mu)\right)\right\}_n \sim_{\lambda} \frac{t_\alpha}{\text{\rm g}_{\alpha}^{(2)}}.
    \end{equation*}
\end{prop}

\begin{prova}
    From the proof of Proposition \ref{prop:negative-preconditioner-1} we know that $\{\tau\big(T_{\mi{n}}(\text{g}_{\alpha}^{(2)})\big)\}_n\sim_{\text{GLT}} \text{g}_{\alpha}^{(2)}$. By {\bf GLT 3}, {\bf GLT 1} and
    Equation \eqref{eq:diagonal-part-precond-distr}, it holds $\big\{\tau\big(T_{\mi{n}}(\text{g}_{\alpha}^{(2)})\big)+ h^\alpha \mu_h I_{n^2}\big\}_n\sim_{\text{GLT},\sigma} \text{g}_{\alpha}^{(2)}$. The conclusion follows as in Proposition \ref{prop:negative-preconditioner-1}. If $\mu(x,y)$ is real-valued, the reasoning is the same, combined with the fact that the diagonal matrix $h^\alpha \mu_h I_{n^2}$ is Hermitian.
\end{prova}

To summarise, both preconditioned matrix sequences cannot cluster at 1 in the sense of the  singular values, since the cluster is at the range of $t_\alpha / \text{g}_{\alpha}^{(2)}$. The same conclusion is drawn for the eigenvalues if $\mu(x,y)$ is real-valued. In light of these negative results, the main part of our analysis will be focused on the second pair of preconditioners, defined simply as
\begin{equation*} 
    \tau(\hat{B}_{\mi{n}}),  \qquad
    \tau(\hat{B}_{\mi{n}})+h^\alpha \mu_h I_{n^2}.
\end{equation*}
where $\mu_h$ is again the arithmetic average of the evaluations of $\mu(x,y)$ over the grid points. We will prove that, under mild assumptions, these preconditioners guarantee the cluster a 1, both in the sense of the eigenvalues and of the singular values.

\subsection{Singular value and spectral analysis} \label{ssec:spectral}

In this subsection we study the singular value and eigenvalue distribution of the latter preconditioning proposals. In Subsection \ref{sssec:bounded-distr} we analyze the case where the wave number $\mu(x,y)$ is essentially bounded. In this case, the distributions follow in a straightforward way from the GLT nature of the involved matrix sequences. Then in Subsection \ref{sssec:unbounded-distr} we progress to the more general case of an unbounded function $\mu(x,y)$. In this case, the eigenvalue distribution is not automatically implied by the GLT axioms on account of the fact that the matrix sequences are not Hermitian.

\subsubsection{The case of an essentially bounded wave function} \label{sssec:bounded-distr}

\begin{theorem}\label{th-main:original-P matrix-sequence SV}
Let $\mu(x,y)$ be an essentially bounded complex-valued function and let $\{\bar{\mu}_{h(n)}\}_n$ be a bounded sequence of complex numbers.
Then, for $P_{\mi{n}}\in \left\{\tau(\hat{B}_{\mi{n}}),\tau(\hat{B}_{\mi{n}})+h^\alpha \bar{\mu}_{h(n)} I_{n^2} \right\}$,
\begin{itemize}
    \item[\bf p1.] $\{P_{\mi{n}}\}_n \sim_{{\rm GLT},\sigma, \lambda} t_\alpha$;
    \item[\bf p2.] $\left\{P_{\mi{n}}^{-1} \left(\hat{B}_{\mi{n}}+h^\alpha D_{\mi{n}}(\mu)\right)\right\}_n \sim_{{\rm GLT},\sigma} 1$. \\
    If moreover $\{\bar{\mu}_{h(n)}\}_n$ is made of real nonnegative numbers, $\left\{P_{\mi{n}}^{-1} \left(\hat{B}_{\mi{n}}+h^\alpha D_{\mi{n}}(\mu)\right)\right\}_n \sim_{\lambda} 1$.
\end{itemize}
\end{theorem}

\begin{prova}
To prove {\bf p1}, we observe that $\{\delta_{\mi{n}}=\hat{B}_{\mi{n}}-\tau(\hat{B}_{\mi{n}})\}_n \sim_{\rm GLT} 0$, since by \eqref{hankel:distr} the sequence is zero-distributed. Therefore by {\bf GLT 3} $\{\tau(\hat{B}_{\mi{n}})=\hat{B}_{\mi{n}}-\delta_{\mi{n}}\}_n\sim_{\rm GLT} t_\alpha - 0=t_\alpha$. By \eqref{eq:diagonal-part-precond-distr}, $\{h^\alpha \bar{\mu}_{h(n)} I_{n^2}\}_n\sim_{\rm GLT} 0$ and again {\bf GLT 3} implies $\{\tau(\hat{B}_{\mi{n}})+h^\alpha \bar{\mu}_{h(n)}I_{n^2}\}_n\sim_{\rm GLT} t_\alpha$. So by {\bf GLT 1} it holds $\{P_{\mi{n}}\}_n\sim_\sigma t_\alpha$ for $P_{\mi{n}}\in \left\{\tau(\hat{B}_{\mi{n}}),\tau(\hat{B}_{\mi{n}})+h^\alpha \bar{\mu}_h I \right\}$. To handle the eigenvalue distribution, we remark that, since $\hat{B}_{\mi{n}}$ is symmetric positive definite, by similarity the desired result is equivalent to
\[
\{X_n+Z_n\}_n\sim_\lambda 1
\]
with
\begin{align*}
X_n & \;:=\; \tau(\hat{B}_{\mi{n}})^{-1/2}\hat{B}_{\mi{n}}\tau(\hat{B}_{\mi{n}})^{-1/2}, \\
Z_n & \;:=\; \tau(\hat{B}_{\mi{n}})^{-1/2}h^\alpha D_{\mi{n}}(\mu)\tau(\hat{B}_{\mi{n}})^{-1/2}.
\end{align*}
To prove the latter, we rely on Theorem \ref{th:GLT-lambda}. Since $\{X_n\}_n$ is a symmetric positive definite GLT sequence with symbol 1, by the second part of {\bf GLT 1} we have $\{X_n\}_n\sim_\lambda 1$ and assumption 2 is satisfied. Regarding assumption 1, by the results in \cite{DB,tau-essential}  $\{X_n\}_n$ is uniformly bounded in spectral norm and the minimal eigenvalue of $\tau(\hat{B}_{\mi{n}})$ is asymptotic to $h^\alpha$, leading to the fact that {$\{Z_n\}_n$} is also uniformly bounded in spectral norm i.e. its maximal singular value is $O(1)$ (see (\ref{first sv}) for more precise derivations). Finally, considering that the eigenvalues of $\tau(\hat{B}_{\mi{n}})$ are a uniform sampling of $t_\alpha$ (up to relative infinitesimal errors), the $s$-th singular value of $Z_n$ is bounded from above by
\[
r_s(\alpha, \mu):=\frac{h^\alpha \|\mu\|_\infty}{\lambda_s\left(\tau(\hat{B}_{\mi{n}})\right)},
\]
so that 
\begin{eqnarray}\label{first sv}
  \sigma_1(Z_n) & \le & r_1(\alpha, \mu)=O(1), \\ \label{almost all sv small}
  \sigma_s(Z_n) & \le & r_s(\alpha, \mu)=o(1),\ \ \forall s=s(n), \ \ 1=0(s(n)).
\end{eqnarray}
Therefore $\|Z_n\|_1=\sum_{j=1}^{n^2}\sigma_j(Z_n)= o(n^2)$ by direct computation using (\ref{first sv}), (\ref{almost all sv small}), and the ordering
$\sigma_1(Z_n)\ge \sigma_j(Z_n)\ge \cdots \ge \sigma_{n^2}(Z_n)$. Applying Theorem \ref{th:GLT-lambda} we obtain $\{Y_n=X_n+Z_n\}_n\sim_\lambda 1$.

The GLT properties in {\bf p2} follow plainly from Theorem \ref{th:coeff matrix-sequence SV} and the third and fourth parts of {\bf GLT 3} combined with {\bf p1}. The singular value distribution then is an obvious consequence of {\bf GLT 1}.
As for the eigenvalue distribution, it follows in the same way as above from Theorem \ref{th:GLT-lambda}, given that the sequence $\bar{\mu}_{h(n)}$ is nonnegative and bounded.
\end{prova}

\subsubsection{The case of an unbounded wave function} \label{sssec:unbounded-distr}

\begin{theorem}\label{th-main-EXT:original-P matrix-sequence}
Let $\mu_M(x,y)$ in (\ref{mu truncated}) be Riemann-integrable for any $M>0$ and suppose that $\{\mu_M(x,y)\}_M$ converges pointwise to $\mu(x,y)$, with $\mu(x,y)$ unbounded and complex-valued. Let $\{\bar{\mu}_{h(n)}\}_n$ be a sequence of complex numbers such that $h^\alpha\bar{\mu}_{h(n)}=o(1)$.
Then, for $P_{\mi{n}}\in \left\{\tau(\hat{B}_{\mi{n}}),\tau(\hat{B}_{\mi{n}})+h^\alpha \bar{\mu}_{h(n)} I_{n^2} \right\}$,
\begin{itemize}
    \item[\bf q1.] $\{P_{\mi{n}}\}_n\sim_{{\rm GLT}, \sigma} t_\alpha$; 
    \item[\bf q2.] $\left\{P_{\mi{n}}^{-1} \left(\hat{B}_{\mi{n}}+h^\alpha D_{\mi{n}}(\mu)\right)\right\}_n \sim_{{\rm GLT}, \sigma} 1$.
\end{itemize}
If moreover $\mu(x,y)$ and $\{\bar{\mu}_{h(n)}\}_n$ are real-valued, then
\begin{itemize}
    \item[\bf q3.] $\{P_{\mi{n}}\}_n\sim_\lambda t_\alpha$;
    \item[\bf q4.] $\left\{P_{\mi{n}}^{-1}\left( \hat{B}_{\mi{n}}+h^\alpha D_{\mi{n}}(\mu)\right)\right\}_n\sim_\lambda 1$.
\end{itemize}
\end{theorem}

\begin{prova}
From Theorem \ref{th-main:original-P matrix-sequence SV} we know that $\{\tau(\hat{B}_{\mi{n}})\}_n \sim_{\rm GLT} t_\alpha$ and from \eqref{eq:diagonal-part-precond-distr} that $\{h^\alpha \bar{\mu}_{h(n)} I_{n^2}\}_n\sim_{\rm GLT} 0$. By {\bf GLT 3}, we have $\{P_{\mi{n}}\}_n \sim_{\rm GLT} t_\alpha$ for $P_{\mi{n}}\in \big\{\tau(\hat{B}_{\mi{n}}),\tau(\hat{B}_{\mi{n}})+h^\alpha \bar{\mu}_h I_{n^2} \big\}$, and by {\bf GLT 1} $\{P_{\mi{n}}\}_n \sim_{\sigma} t_\alpha$.

The GLT properties in {\bf q2} follow plainly from Theorem \ref{th:coeff matrix-sequence SV-EXT} and the third and fourth parts of {\bf GLT 3} combined with {\bf q1}.

If $\mu$ and $\{\bar{\mu}_{h(n)}\}_n$ are real, all the matrices are Hermitian or similar to Hermitian matrices, where the similarity is still of GLT nature. Therefore statements {\bf q3} and {\bf q4} follow from the second part of {\bf GLT 1}.
\end{prova}

As we have already mentioned, when $\mu(x,y)$ is complex-valued and unbounded we have been unable to obtain a satisfactory generalization of the eigenvalue distribution. In particular, the unbounded character of the inverse of the preconditioning matrix-sequences has been the main obstacle. This matter is therefore left for future work, see in particular \cite{AST} for the specific case of a complex-valued $\mu$ with power singularities.

\section{Numerical experiments}\label{sec:num-exp}

This section is divided in two parts. Subsection \ref{ssec:visualization} contains a visualization of the most relevant theoretical results, while
Subsection \ref{ssec:p-gmres} analyzes the convergence behavior of the preconditioned CG/GMRES applied to the scaled linear systems \eqref{scaled eq}
\begin{equation*}
    \hat{A}_{\mi{n}} \vet{u} := (\hat{B}_{\mi{n}}+h^\alpha D_{\mi{n}}(\mu)) \vet{u} = h^\alpha \vet{v}, \qquad \mi{n}=(n,n),
\end{equation*}
with the preconditioners described in Section \ref{ssec:spectral}
\begin{equation*}
    P_{\mi{n}} = \tau(\hat{B}_{\mi{n}})+h^\alpha \mu_h I_{n^2},
  \end{equation*}
where $\hat{B}_{\mi{n}} := T_{\mi{n}}(t_\alpha)$ and ${\mu}_h$ is chosen as the arithmetic average of the evaluations of $\mu(x,y)$ over the grid points.

\subsection{Numerical tests and visualization of the theory}\label{ssec:visualization}

In the current subsection we report a numerical and visual representation the singular value and eigenvalue analysis made in the previous sections.

\begin{figure}
\centering
  \includegraphics[width=\textwidth]{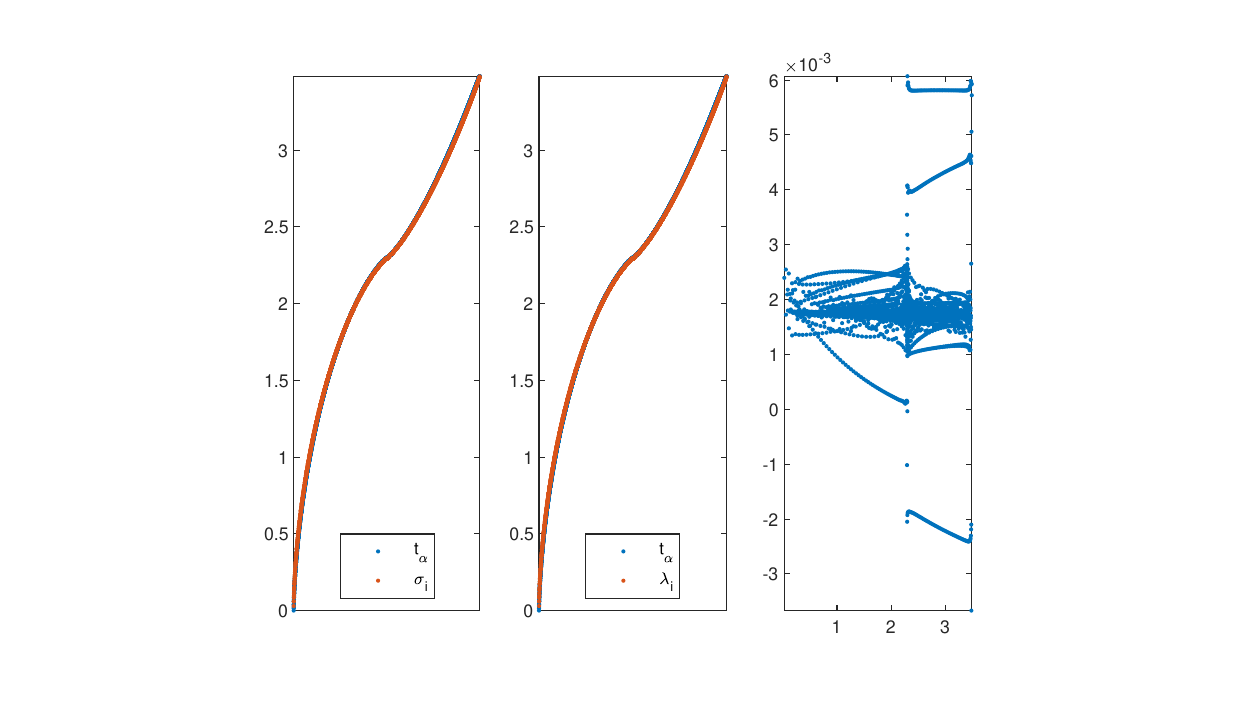}
  \caption{Singular/Eigenvalues of the  matrix $\hat{A}_{\mi{n}}$ for $\mu(x,y)=\exp(\iu(x+4y))$, $\alpha=1.2$ and $n^2=2^{12}$. The first (second) panel reports in blue the singular values (real part of the eigenvalues) and in red the equispaced sampling of $t_\alpha$ in nondecreasing order. The third panel reports the eigenvalues in the complex plane.}
\label{Fig_flagmu11_12}
\end{figure}

\begin{figure}
\centering
  \includegraphics[width=\textwidth]{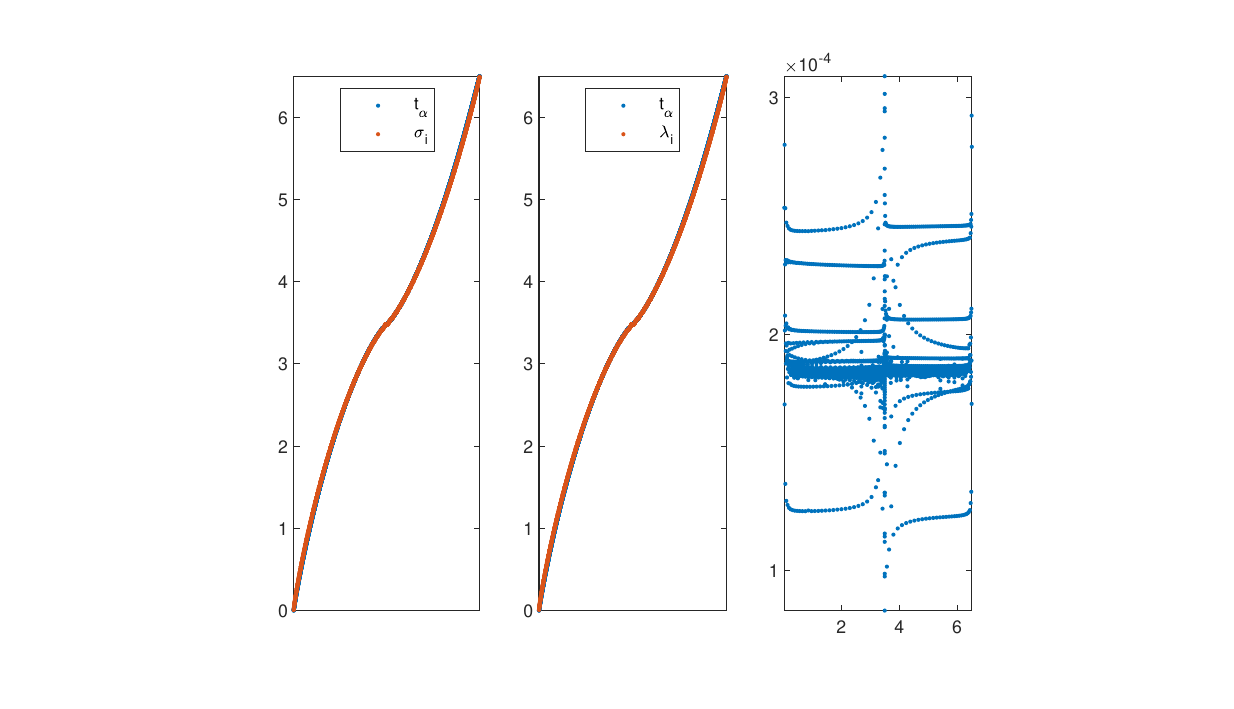}
  \caption{Singular/Eigenvalues of the  matrix $\hat{A}_{\mi{n}}$ for $\mu(x,y)=-2+\exp(\iu(3x+2y))$, $\alpha=1.8$ and $n^2=2^{12}$.  The first (second) panel reports in blue the singular values (real part of the eigenvalues) and in red the equispaced sampling of $t_\alpha$ in nondecreasing order. The third panel reports the eigenvalues in the complex plane.}
\label{Fig_flagmu12_18}
\end{figure}

In Figures \ref{Fig_flagmu11_12}--\ref{Fig_flagmu12_18}, we plot the singular values and eigenvalues of the matrix $\hat{A}_{\mi{n}}$ with size $\nu(\mi{n})=n^2=2^{12}$ and $\alpha\in \{1.2, 1.4, 1.6, 1.8\}$, displaying the results contained in Theorems \ref{th:coeff matrix-sequence SV} and \ref{th:coeff matrix-sequence SV-EXT}. The clustering at zero of the imaginary part of the eigenvalues of $\hat{A}_{\mi{n}}$ and the relation $\{\hat{A}_{\mi{n}}\}\sim_{\sigma,\lambda} t_\alpha$ are already clearly visible for a moderate matrix size like $2^{12}$. It is remarkable that no outliers are present, since the imaginary parts are always negligible and the graph of the equispaced sampling of $t_\alpha$ and of the real parts of the eigenvalues superpose completely. The same happens for the singular values.

Tables \ref{tab:as_flag11svd}--\ref{tab:as_flag12svd} show the number of singular values outliers, with respect to neighbourhoods of radius $\epsilon = 0.1, \, 0.01$, for the preconditioned matrix sequence, for $\alpha\in \{1.2, 1.4, 1.6, 1.8\}$ and two examples of the function $\mu$. The clustering at $1$ is evident and in accordance with Theorem \ref{th-main-EXT:original-P matrix-sequence}, item {\bf q2}.

\begin{figure}
\centering
  \includegraphics[width=\textwidth]{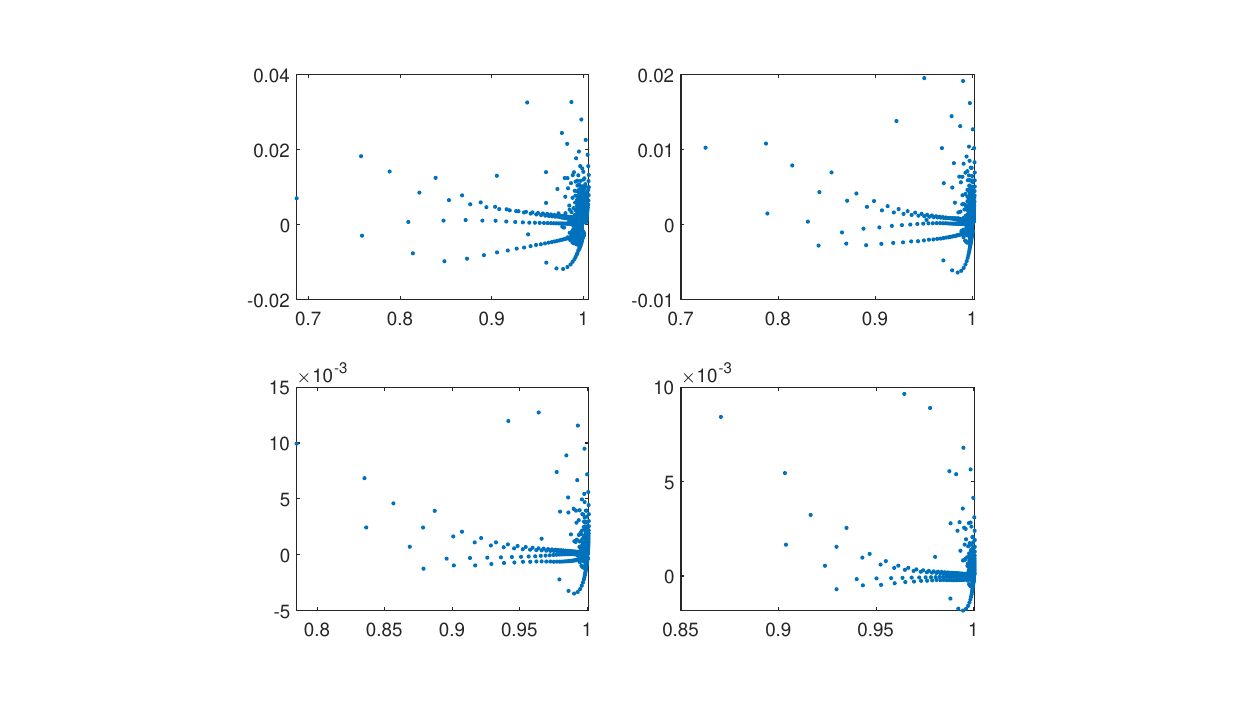}
  \caption{Eigenvalues of the preconditioned matrix of size $n^2=2^{12}$ for $\mu(x,y)=\exp(\iu(x+4y))$ and $\alpha=\{1.2,1.4,1.6,1.8\}$, respectively.}
\label{Fig_Preig_flagmu11}
\end{figure}

\begin{figure}
\centering
  \includegraphics[width=\textwidth]{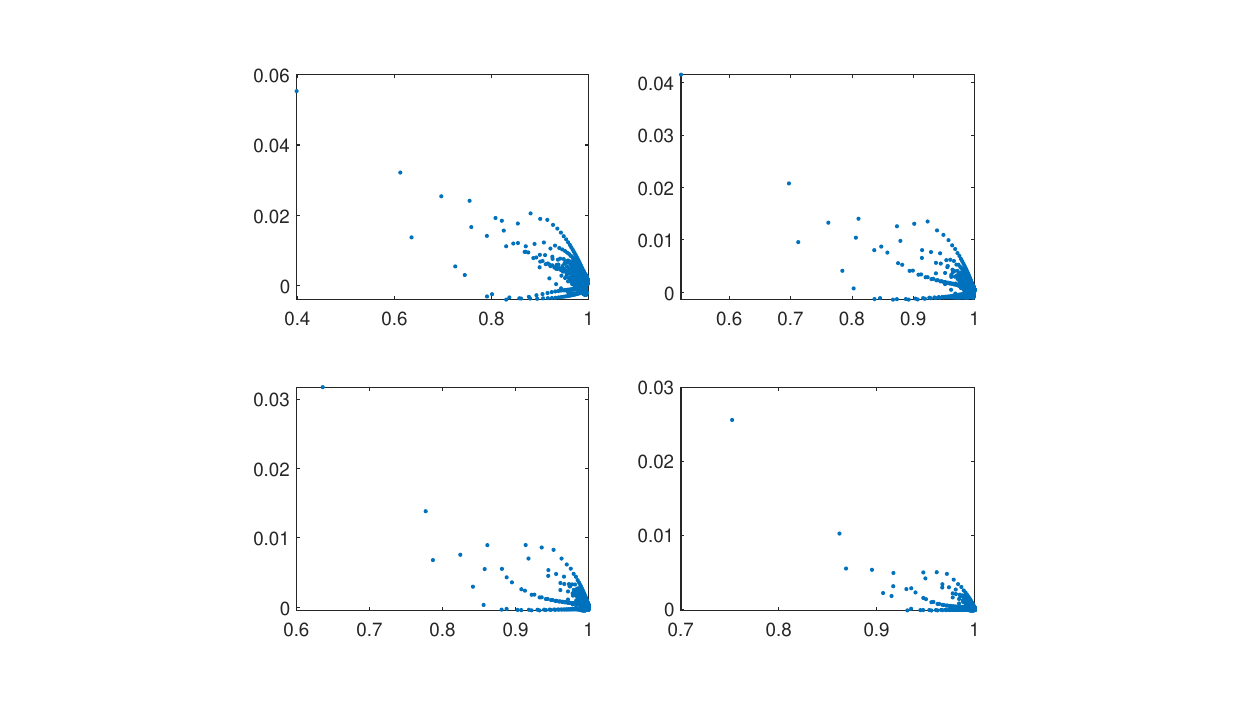}
  \caption{Eigenvalues of the preconditioned matrix of size $n^2=2^{12}$ for $\mu(x,y)=-2+\exp(\iu(3x+2y))$, and $\alpha=\{1.2,1.4,1.6,1.8\}$, respectively.}
\label{Fig_Preig_flagmu12}
\end{figure}

The same kind of evidence for the outlying eigenvalues is given in Tables \ref{tab:as_flag13}--\ref{tab:as_flag14}. The spectral clustering at 1 emerges from the numerical data and it is consistent with the theoretical analysis in {Theorems \ref{th-main:original-P matrix-sequence SV}--\ref{th-main-EXT:original-P matrix-sequence}}. Figures \ref{Fig_Preig_flagmu11}--\ref{Fig_Preig_flagmu12} provide a further visualization of the eigenvalue clustering at 1, showing the whole spectra in the complex plane.

\begin{figure}
\centering
  \includegraphics[width=\textwidth]{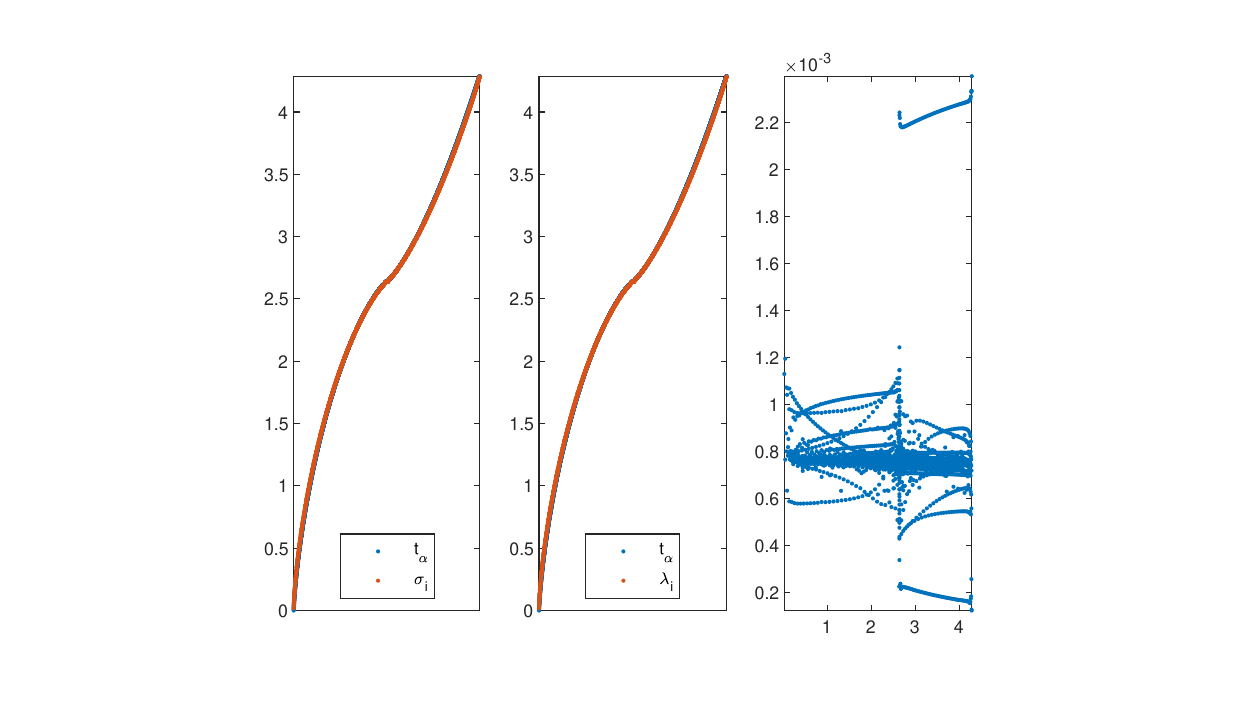}
  \caption{Singular/Eigenvalues of the  matrix $\hat{A}_{\mi{n}}$ for $\mu(x,y)=\exp(\iu(x+4y))$, $\alpha=1.4$ and $n^2=2^{12}$.  The first (second) panel reports in blue the singular values (real part of the eigenvalues) and in red the equispaced sampling of $t_\alpha$ in nondecreasing order.   The third panel reports the eigenvalues in the complex plane.}
\end{figure}
\begin{figure}
\centering
  \includegraphics[width=\textwidth]{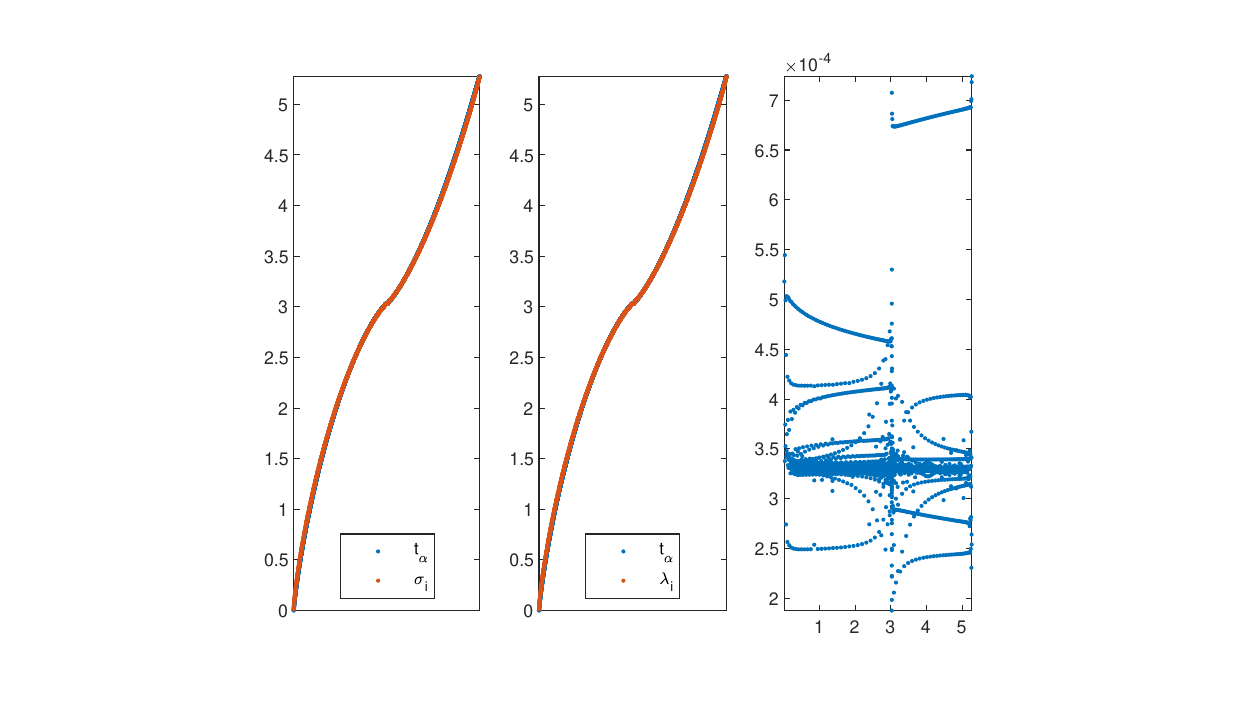}
  \caption{Singular/Eigenvalues of the  matrix $\hat{A}_{\mi{n}}$ for $\mu(x,y)=\exp(\iu(x+4y))$, $\alpha=1.6$ and $n^2=2^{12}$.  The first (second) panel reports in blue the singular values (real part of the eigenvalues) and in red the equispaced sampling of $t_\alpha$ in nondecreasing order. The third panel reports the eigenvalues in the complex plane.}
\end{figure}
\begin{figure}
\centering
  \includegraphics[width=\textwidth]{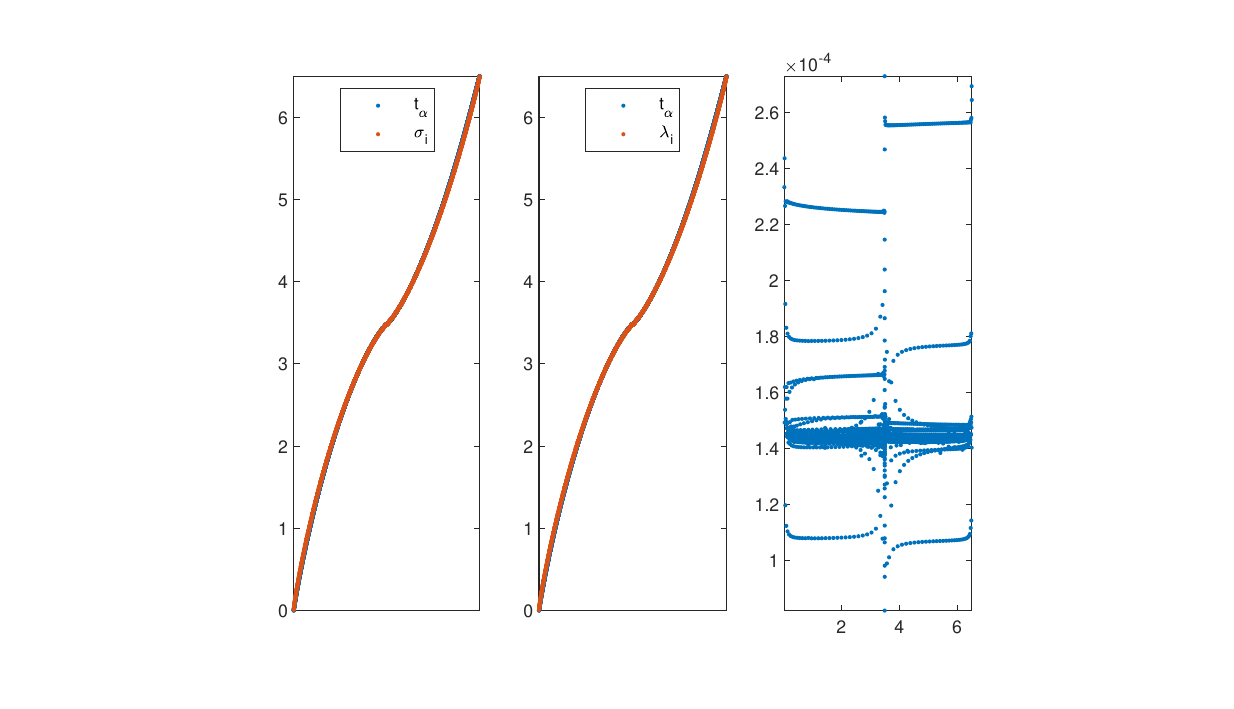}
  \caption{Singular/Eigenvalues of the  matrix $\hat{A}_{\mi{n}}$ for $\mu(x,y)=\exp(\iu(x+4y))$, $\alpha=1.8$ and $n^2=2^{12}$.  The first (second) panel reports in blue the singular values (real part of the eigenvalues) and in red the equispaced sampling of $t_\alpha$ in nondecreasing order. The third panel reports the eigenvalues in the complex plane.}
\end{figure}
\begin{figure}
\centering
  \includegraphics[width=\textwidth]{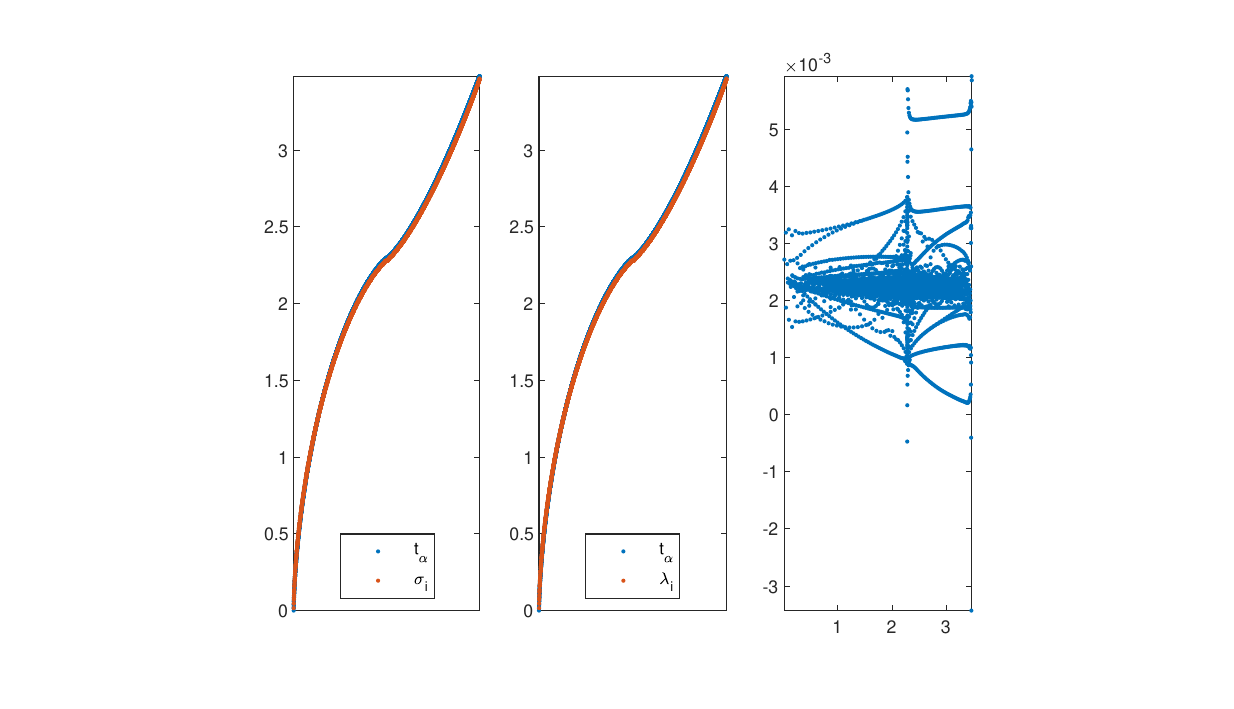}
  \caption{Singular/Eigenvalues of the  matrix $\hat{A}_{\mi{n}}$ for $\mu(x,y)=-2+\exp(\iu(3x+2y))$, $\alpha=1.2$ and $n^2=2^{12}$.  The first (second) panel reports in blue the singular values (real part of the eigenvalues) and in red the equispaced sampling of $t_\alpha$ in nondecreasing order. The third panel reports the eigenvalues in the complex plane.}
\end{figure}
\begin{figure}
\centering
  \includegraphics[width=\textwidth]{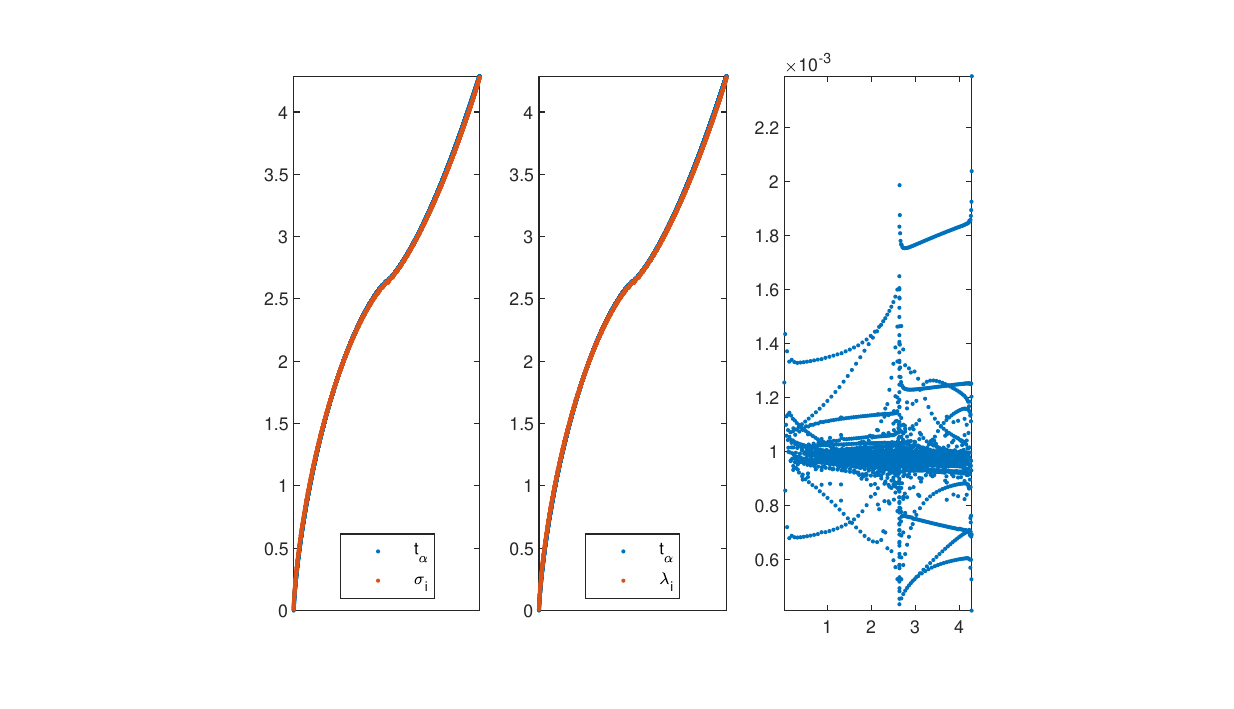}
  \caption{Singular/Eigenvalues of the  matrix $\hat{A}_{\mi{n}}$ for $\mu(x,y)=-2+\exp(\iu(3x+2y))$, $\alpha=1.4$ and $n^2=2^{12}$.  The first (second) panel reports in blue the singular values (real part of the eigenvalues) and in red the equispaced sampling of $t_\alpha$ in nondecreasing order. The third panel reports the eigenvalues in the complex plane.}
\end{figure}
\begin{figure}
\centering
  \includegraphics[width=\textwidth]{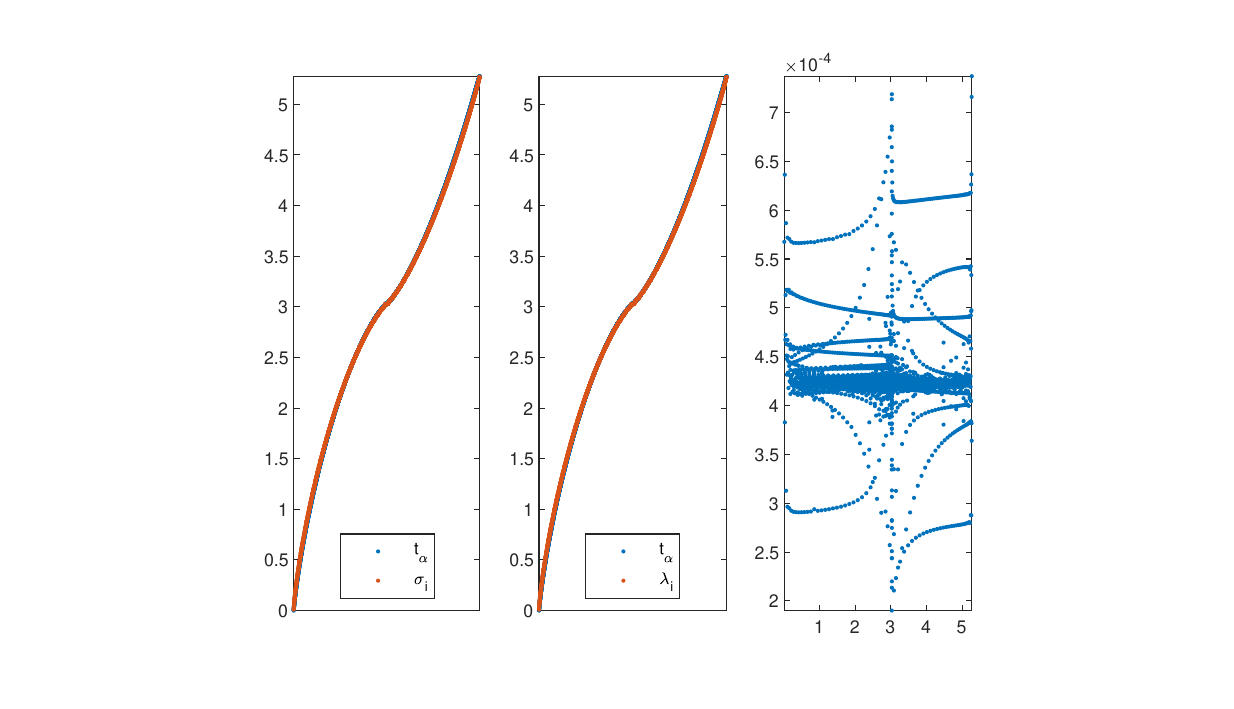}
  \caption{Singular/Eigenvalues of the matrix $\hat{A}_{\mi{n}}$ for $\mu(x,y)=-2+\exp(\iu(3x+2y))$, $\alpha=1.6$ and $n^2=2^{12}$.  The first (second) panel reports in blue the singular values (real part of the eigenvalues) and in red the equispaced sampling of $t_\alpha$ in nondecreasing order. The third panel reports the eigenvalues in the complex plane.}
\end{figure}
%
\begin{table}
  \centering
  \footnotesize
  \begin{tabular}{|c |c|c |c|c |}
\hline
\multicolumn {5}{|c|}{$\mu(x,y)=\exp(\iu(x+4y))$} \\[5pt]
$n$ & $N_o(0.1)$ & $\rm{Percentage}$ & $N_o(0.01)$ & $\rm{Percentage}$  \\
\hline
     &  \multicolumn {4}{|c|}{$\alpha=1.2$}\\
$2^4$ & 5  & 1.95\% & 130 & 50.8\% \\
$2^5$ & 11 & 1.07\% & 137 & 13.4\% \\
$2^6$ & 21 & 0.51\% & 191 & 4.66\% \\
     \hline
     &  \multicolumn {4}{|c|}{$\alpha=1.4$}\\
$2^4$ & 4  & 1.56\% & 51 & 20.0\% \\
$2^5$ & 9  & 0.88\% & 77 & 7.52\% \\
$2^6$ & 19 & 0.46\% & 136 & 3.32\% \\
     \hline
     &  \multicolumn {4}{|c|}{$\alpha=1.6$}\\
$2^4$ & 3 & 1.17\% & 31 & 12.1\% \\
$2^5$ & 6 & 0.58\% & 57 & 5.57\% \\
$2^6$ & 13 & 0.32\% & 114 & 2.78\% \\	
     \hline
     &  \multicolumn {4}{|c|}{$\alpha=1.8$}\\
$2^4$ & 1 & 0.39\% & 20 & 7.81\% \\
$2^5$ & 2 & 0.19\% & 40 & 3.91\% \\
$2^6$ & 5 & 0.12\% & 82 & 2.00\% \\
     \hline
  \end{tabular}
  \caption{Number of outlying singular values $N_o(\epsilon)$ with respect to a neighborhood of $1$ of radius $\epsilon=0.1$ or $\epsilon=0.01$ and related percentage for increasing dimension $n^2$ - $\Omega=[0,1]^2$.} \label{tab:as_flag11svd}
\end{table}
%
\begin{table}
  \centering
  \footnotesize
  \begin{tabular}{|c |c|c |c|c |}
\hline
\multicolumn {5}{|c|}{$\mu(x,y)=-2+\exp(\iu(3x+2y))$} \\
$n$ & $N_o(0.1)$ & $\rm{Percentage}$ & $N_o(0.01)$ & $\rm{Percentage}$  \\
\hline
     &  \multicolumn {4}{|c|}{$\alpha=1.2$}\\
$2^4$ & 20 & 7.81\% & 256 & 100\% \\
$2^5$ & 25 & 2.44\% & 932 & 91.0\% \\
$2^6$ & 36 & 0.88\% & 1065 & 26.0\% \\
     \hline
     &  \multicolumn {4}{|c|}{$\alpha=1.4$}\\
$2^4$ & 10 & 3.91\% & 243 & 94.9\% \\
$2^5$ & 15 & 1.46\% & 297 & 29.0\% \\
$2^6$ & 25 & 6.10\% & 304 & 7.42\% \\
     \hline
     &  \multicolumn {4}{|c|}{$\alpha=1.6$}\\
$2^4$ & 5 & 1.95\% & 122 & 47.6\% \\
$2^5$ & 8 & 0.78\% & 114 & 11.1\% \\
$2^6$ & 16 & 0.39\% & 160 & 3.91\% \\
     \hline
     &  \multicolumn {4}{|c|}{$\alpha=1.8$}\\
$2^4$ & 3 & 1.17\% & 47 & 18.3\% \\
$2^5$ & 4 & 0.39\% & 60 & 5.86\% \\
$2^6$ & 7 & 0.17\% & 99 & 2.42\% \\
     \hline
  \end{tabular}
  \caption{Number of outlying singular values $N_o(\epsilon)$ with respect to a neighborhood of $1$ of radius $\epsilon=0.1$ or $\epsilon=0.01$ and related percentage for increasing dimension $n^2$ - $\Omega=[0,1]^2$.} \label{tab:as_flag12svd}
\end{table}

\begin{table}
  \centering
  \footnotesize
  \begin{tabular}{|c |c|c |c|c |}
\hline
\multicolumn {5}{|c|}{$\mu(x,y)=\exp(\iu(x+4y))$} \\[5pt]
$n$ & $N_o(0.1)$ & $\rm{Percentage}$ & $N_o(0.01)$ & $\rm{Percentage}$  \\
\hline
     &  \multicolumn {4}{|c|}{$\alpha=1.2$}\\
$2^4$ & 5 & 1.95\% & 247  & 96.5\% \\
$2^5$ & 10 & 0.98\% & 165 & 16.1\% \\
$2^6$ & 19 & 0.46\% & 211 & 5.15\% \\
     \hline
     &  \multicolumn {4}{|c|}{$\alpha=1.4$}\\
$2^4$ & 4 & 1.56\% & 71 & 27.7\% \\
$2^5$ & 7 & 0.68\% & 83 & 8.10\% \\
$2^6$ & 16 & 0.390\% & 142 & 3.47\% \\
     \hline
     &  \multicolumn {4}{|c|}{$\alpha=1.6$}\\
$2^4$ & 3 & 1.17\% & 34 & 13.3\% \\
$2^5$ & 4 & 0.39\% & 56 & 5.47\% \\
$2^6$ & 9 & 0.22\% & 109 & 2.66\% \\
     \hline
     &  \multicolumn {4}{|c|}{$\alpha=1.8$}\\
$2^4$ & 0 & 0\% & 21 & 8.20\% \\
$2^5$ & 1 & 0.09\% & 38 & 3.71\% \\
$2^6$ & 1 & 0.02\% & 75 & 1.83\% \\
     \hline
  \end{tabular}
  \caption{Number of outlying eigenvalues $N_o(\epsilon)$ with respect to a neighborhood of $1$ of radius $\epsilon=0.1$ or $\epsilon=0.01$ and related percentage for increasing dimension $n^2$ - $\Omega=[0,1]^2$.} \label{tab:as_flag13}
\end{table}
%
\begin{table}
  \centering
  \footnotesize
  \begin{tabular}{|c |c|c |c|c |}
\hline
\multicolumn {5}{|c|}{$\mu(x,y)=-2+\exp(\iu(3x+2y))$} \\[5pt]
$n$ & $N_o(0.1)$ & $\rm{Percentage}$ & $N_o(0.01)$ & $\rm{Percentage}$  \\
\hline
     &  \multicolumn {4}{|c|}{$\alpha=1.2$}\\
$2^4$ & 20 & 7.81\% & 256 & 100\% \\
$2^5$ & 25 & 2.44\% & 968 & 94.5\% \\
$2^6$ & 38 & 0.93\% & 1078 & 26.3\% \\
     \hline
     &  \multicolumn {4}{|c|}{$\alpha=1.4$}\\
$2^4$ & 10 & 3.91\% & 249 & 97.2\% \\
$2^5$ & 15 & 1.46\% & 304 & 29.7\% \\
$2^6$ & 23 & 0.56\% & 323 & 7.88\% \\
     \hline
     &  \multicolumn {4}{|c|}{$\alpha=1.6$}\\
$2^4$ & 6 & 2.34\% & 126 & 49.2\% \\
$2^5$ & 8 & 0.78\% & 121 & 11.8\% \\
$2^6$ & 13 & 0.32\% & 173 & 4.22\% \\
     \hline
     &  \multicolumn {4}{|c|}{$\alpha=1.8$}\\
$2^4$ & 3 & 1.17\% & 47 & 18.3\% \\
$2^5$ & 3 & 0.29\% & 64 & 6.25\% \\
$2^6$ & 4 & 0.09\% & 99 & 2.42\% \\
     \hline
  \end{tabular}
  \caption{Number of outlying eigenvalues $N_o(\epsilon)$ with respect to a neighborhood of $1$ of radius $\epsilon=0.1$ or $\epsilon=0.01$ and related percentage for increasing dimension $n^2$ - $\Omega=[0,1]^2$.} \label{tab:as_flag14}
\end{table}
%

\subsection{Preconditioned CG and GMRES convergence}\label{ssec:p-gmres}

In the present subsection, we provide convergence results when applying the considered preconditioning strategy for two Krylov methods to the non-Hermitian linear systems arising from the numerical approximation of our fractional Helmholtz equations on the square domain $\Omega=[L,R]^2=[-1,1]^2$. Here $h = (R-L)/(n+1)$ is the step size, $x_i = L+ih, \ i=0,\ldots,n+1$, $y_j = L+ jh, \ j=0,\ldots,n+1$, and the domain $\Omega=[-1,1]^2$ is chosen as in the numerical results in \cite{tau-prec}.


As expected from {Theorem \ref{th-main:original-P matrix-sequence SV}, item {\bf p2}, and Theorem \ref{th-main-EXT:original-P matrix-sequence}, item {\bf q4}}, when $\mu(x,y)$ is real-valued a standard preconditioned CG (PCG) method is sufficient for obtaining a fast convergence, with a total cost of $O\big(n^2\log(n)\big)$ arithmetic operations. In fact, under mild assumptions, the resulting matrices are real symmetric and positive definite and the associated preconditioned matrix-sequences show a spectral clustering at $1$. The fast convergence is well confirmed in Table \ref{tab:pcg_flag4567}, Table \ref{tab:pcg_flag89} and in Table \ref{tab:pcggmres_flag10}. Of course, when the infinity norm of $\mu_M(x,y)=M \cos(x)\cos(y)$ grows ($M=1,10,100$), we observe a slight increase in the number of the iterations, just because the radius of the cluster at 1 is influenced by $M$ as the matrix size is fixed. In any case, we can state that the resulting PCG is optimal and robust with respect to the problem and discretization parameters $\alpha, M, h$.
In Table \ref{tab:gmres_flag89} and in Table \ref{tab:pcggmres_flag10} we report the numerical results when a preconditioned GMRES is employed with the same preconditioner and the quality of the convergence is of the same type as for the PCG.

Furthermore, in Table \ref{tab:gmres_flag11} and in Table \ref{tab:gmres_flag12}, we take into consideration two examples of bounded complex-valued functions $\mu(x,y)$. In accordance with the theoretical analysis and with the eigenvalue clustering at $1$ observed in Subsection \ref{ssec:visualization}, a very good numerical behavior is reported. Also for these two examples the preconditioned GMRES is optimal and robust with respect to the parameters $\alpha$ and $h$.

Finally, we stress that the results are remarkable, {especially} taking into account the asymptotic theoretical barriers proved in \cite{prec-neg,nega-gen} concerning the matrix algebra preconditioning for linear systems with multilevel Toeplitz coefficient matrices.

\begin{table}
  \centering
  \footnotesize
  \begin{tabular}{|c |c|c |c|c |c|c |c|c|}
  \hline
\multicolumn {9}{|c|}{$\mu(x,y)=\cos(x)\cos(y)$} \\
\hline
     &  \multicolumn {2}{|c|}{$\alpha=1.2$} &  \multicolumn {2}{|c|}{$\alpha=1.4$} &  \multicolumn {2}{|c|}{$\alpha=1.6$} &  \multicolumn {2}{|c|}{$\alpha=1.8$}\\
$n$ & - & $P_{\mi{n}}$ & - & $P_{\mi{n}}$ & - & $P_{\mi{n}}$ & - & $P_{\mi{n}}$   \\
   \hline
$2^4$ & 23 & 8 & 25 & 8 & 28 & 7 & 30 & 7 \\
$2^5$ & 36 & 9 & 44 & 9 & 51 & 8 & 60 & 7 \\
$2^6$ & 56 & 10 & 72 & 10 & 91 & 9 & 113 & 8 \\
$2^7$ & 85 & 11 & 115 & 10 & 159 & 10 & 212 & 8 \\
$2^8$ & 131 & 12 & 190 & 11 & 276 & 10 & 397 & 9 \\
$2^9$ & 200 & 12 & 314 & 12 & 489 & 11 & $>$500 & 9 \\
   \hline
\multicolumn {9}{|c|}{$\mu(x,y)=\cos(2x)\cos(2y)$} \\
\hline
     &  \multicolumn {2}{|c|}{$\alpha=1.2$} &  \multicolumn {2}{|c|}{$\alpha=1.4$} &  \multicolumn {2}{|c|}{$\alpha=1.6$} &  \multicolumn {2}{|c|}{$\alpha=1.8$}\\
$n$ & - & $P_{\mi{n}}$ & - & $P_{\mi{n}}$ & - & $P_{\mi{n}}$ & - & $P_{\mi{n}}$   \\
   \hline
$2^4$ & 23 & 9 & 25 & 8 & 28 & 8 & 30 & 7 \\
$2^5$ & 37 & 9 & 44 & 9 & 52 & 8 & 60 & 7 \\
$2^6$ & 56 & 10 & 73 & 9 & 91 & 9 & 114 & 8 \\
$2^7$ & 86 & 11 & 116 & 10 & 160 & 9 & 213 & 8 \\
$2^8$ & 132 & 12 & 191 & 11 & 279 & 10 & 399 & 9 \\
$2^9$ & 202 & 12 & 315 & 12 & 492 & 11 & $>$500 & 9 \\
   \hline
\multicolumn {9}{|c|}{$\mu(x,y)=\cos(4x)\cos(4y)$} \\
\hline
     &  \multicolumn {2}{|c|}{$\alpha=1.2$} &  \multicolumn {2}{|c|}{$\alpha=1.4$} &  \multicolumn {2}{|c|}{$\alpha=1.6$} &  \multicolumn {2}{|c|}{$\alpha=1.8$}\\
$n$ & - & $P_{\mi{n}}$ & - & $P_{\mi{n}}$ & - & $P_{\mi{n}}$ & - & $P_{\mi{n}}$   \\
   \hline
$2^4$ & 23 & 9 & 26 & 8 & 28 & 7 & 30 & 7 \\
$2^5$ & 38 & 9 & 43 & 9 & 52 & 8 & 60 & 7 \\
$2^6$ & 59 & 10 & 74 & 9 & 92 & 9 & 114 & 8 \\
$2^7$ & 90 & 11 & 121 & 10 & 161 & 9 & 214 & 8 \\
$2^8$ & 138 & 11 & 199 & 11 & 282 & 10 & 401 & 9 \\
$2^9$ & 211 & 12 & 324 & 11 & 495 & 11 & $>$500 & 9 \\
\hline
\multicolumn {9}{|c|}{$\mu(x,y)=\cos(10x)\cos(10y)$} \\
\hline
     &  \multicolumn {2}{|c|}{$\alpha=1.2$} &  \multicolumn {2}{|c|}{$\alpha=1.4$} &  \multicolumn {2}{|c|}{$\alpha=1.6$} &  \multicolumn {2}{|c|}{$\alpha=1.8$}\\
$n$ & - & $P_{\mi{n}}$ & - & $P_{\mi{n}}$ & - & $P_{\mi{n}}$ & - & $P_{\mi{n}}$   \\
   \hline
$2^4$ & 23 & 8 & 26 & 8 & 28 & 7 & 30 & 7 \\
$2^5$ & 38 & 9 & 45 & 8 & 52 & 8 & 60 & 7 \\
$2^6$ & 59 & 10 & 74 & 9 & 92 & 8 & 114 & 8 \\
$2^7$ & 90 & 10 & 121 & 10 & 162 & 9 & 214 & 8 \\
$2^8$ & 138 & 11 & 199 & 11 & 283 & 10 & 401 & 9 \\
$2^9$ & 210 & 12 & 325 & 11 & 496 & 10 & $>$500 & 9 \\
\hline
  \end{tabular}
  \caption{Number of preconditioned PCG iterations to solve the linear system for increasing dimension $n^2$  till $tol = 1.e-11$ - $\Omega=[-1,1]^2$.} \label{tab:pcg_flag4567}
\end{table}
\begin{table}
  \centering
  \footnotesize
  \begin{tabular}{|c |c|c |c|c |c|c |c|c|}
  \hline
\multicolumn {9}{|c|}{$\mu(x,y)=M\cos(x)\cos(y), \; M=10$} \\
\hline
     &  \multicolumn {2}{|c|}{$\alpha=1.2$} &  \multicolumn {2}{|c|}{$\alpha=1.4$} &  \multicolumn {2}{|c|}{$\alpha=1.6$} &  \multicolumn {2}{|c|}{$\alpha=1.8$}\\
$n$ & - & $P_{\mi{n}}$ & - & $P_{\mi{n}}$ & - & $P_{\mi{n}}$ & - & $P_{\mi{n}}$   \\
   \hline
$2^4$ & 19 & 12 & 22 & 12 & 25 & 11 & 28 & 10 \\
$2^5$ & 29 & 13 & 37 & 12 & 46 & 11 & 56 & 10 \\
$2^6$ & 45 & 14 & 61 & 13 & 82 & 12 & 105 & 11 \\
$2^7$ & 69 & 16 & 100 & 14 & 144 & 12 & 198 & 11 \\
$2^8$ & 106 & 17 & 164 & 15 & 252 & 13 & 371 & 11 \\
$2^9$ & 162 & 17 & 269 & 16 & 442 & 14 & $>$500 & 12 \\
   \hline
\multicolumn {9}{|c|}{$\mu(x,y)=M\cos(x)\cos(y), \; M=100$} \\
\hline
     &  \multicolumn {2}{|c|}{$\alpha=1.2$} &  \multicolumn {2}{|c|}{$\alpha=1.4$} &  \multicolumn {2}{|c|}{$\alpha=1.6$} &  \multicolumn {2}{|c|}{$\alpha=1.8$}\\
$n$ & - & $P_{\mi{n}}$ & - & $P_{\mi{n}}$ & - & $P_{\mi{n}}$ & - & $P_{\mi{n}}$   \\
\hline
$2^4$ & 13 & 19 & 13 & 19 & 16 & 19 & 20 & 18 \\
$2^5$ & 15 & 26 & 20 & 24 & 28 & 22 & 38 & 20 \\
$2^6$ & 22 & 30 & 33 & 28 & 49 & 24 & 71 & 20 \\
$2^7$ & 33 & 33 & 53 & 30 & 86 & 26 & 133 & 22 \\
$2^8$ & 51 & 36 & 88 & 31 & 151 & 27 & 251 & 22 \\
$2^9$ & 78 & 39 & 144 & 33 & 266 & 29 & 471 & 24 \\
\hline
  \end{tabular}
  \caption{Number of preconditioned PCG iterations to solve the linear system for increasing dimension $n^2$  till $tol =1.e-11$ - $\Omega=[-1,1]^2$.} \label{tab:pcg_flag89}
\end{table}
\begin{table}
  \centering
  \footnotesize
  \begin{tabular}{|c |c|c |c|c |c|c |c|c|}
  \hline
\multicolumn {9}{|c|}{$\mu(x,y)=M\cos(x)\cos(y), \; M=10$} \\
\hline
     &  \multicolumn {2}{|c|}{$\alpha=1.2$} &  \multicolumn {2}{|c|}{$\alpha=1.4$} &  \multicolumn {2}{|c|}{$\alpha=1.6$} &  \multicolumn {2}{|c|}{$\alpha=1.8$}\\
$n$ & - & $P_{\mi{n}}$ & - & $P_{\mi{n}}$ & - & $P_{\mi{n}}$ & - & $P_{\mi{n}}$   \\
   \hline
$2^4$ & 19 & 12 & 22 & 11 & 25 & 10 & 28 & 9 \\
$2^5$ & 29 & 12 & 37 & 11 & 46 & 10 & 56 & 9 \\
$2^6$ & 45 & 13 & 60 & 12 & 81 & 11 & 105 & 9 \\
$2^7$ & 68 & 14 & 99 & 12 & 142 & 11 & 195 & 10 \\
$2^8$ & 104 & 14 & 161 & 13 & 248 & 11 & 365 & 10 \\
$2^9$ & 158 & 15 & 263 & 13 & 433 & 11 & $>$500 & 10 \\
   \hline
\multicolumn {9}{|c|}{$\mu(x,y)=M\cos(x)\cos(y), \; M=100$} \\
\hline
     &  \multicolumn {2}{|c|}{$\alpha=1.2$} &  \multicolumn {2}{|c|}{$\alpha=1.4$} &  \multicolumn {2}{|c|}{$\alpha=1.6$} &  \multicolumn {2}{|c|}{$\alpha=1.8$}\\
$n$ & - & $P_{\mi{n}}$ & - & $P_{\mi{n}}$ & - & $P_{\mi{n}}$ & - & $P_{\mi{n}}$   \\
   \hline
$2^4$ & 13 & 17 & 13 & 17 & 16 & 17 & 20 & 16 \\
$2^5$ & 15 & 22 & 20 & 21 & 28 & 19 & 38 & 16 \\
$2^6$ & 22 & 26 & 33 & 23 & 49 & 20 & 71 & 17 \\
$2^7$ & 33 & 27 & 53 & 24 & 85 & 20 & 132 & 18 \\
$2^8$ & 51 & 29 & 87 & 25 & 148 & 21 & 246 & 18 \\
$2^9$ & 77 & 30 & 141 & 26 & 259 & 22 & 462 & 18 \\
   \hline
  \end{tabular}
  \caption{Number of preconditioned GMRES iterations to solve the linear system for increasing dimension $n^2$  till $tol =1.e-11$ - $\Omega=[-1,1]^2$.} \label{tab:gmres_flag89}
\end{table}
\begin{table}
  \centering
  \footnotesize
  \begin{tabular}{|c |c|c |c|c |c|c |c|c|}
\hline
\multicolumn {9}{|c|}{$\mu(x,y)=\cos(x)\cos(4y)(1-\exp(x+y))$ - PCG} \\
\hline
     &  \multicolumn {2}{|c|}{$\alpha=1.2$} &  \multicolumn {2}{|c|}{$\alpha=1.4$} &  \multicolumn {2}{|c|}{$\alpha=1.6$} &  \multicolumn {2}{|c|}{$\alpha=1.8$}\\
$n$ & - & $P_{\mi{n}}$ & - & $P_{\mi{n}}$ & - & $P_{\mi{n}}$ & - & $P_{\mi{n}}$   \\
   \hline
$2^4$ & 41 & 16 & 46 & 14 & 52 & 12 & 60 & 11 \\
$2^5$ & 62 & 16 & 74 & 14 & 92 & 12 & 111 & 11 \\
$2^6$ & 95 & 16 & 122 & 14 & 160 & 13 & 207 & 11 \\
$2^7$ & 145 & 17 & 199 & 15 & 280 & 13 & 389 & 11 \\
$2^8$ & 221 & 17 & 326 & 15 & 491 & 13 & $>$500 & 11 \\
$2^9$ & 338 & 18 & $>$500 & 16 & $>$500 & 14 & $>$500 & 12 \\
\hline
\multicolumn {9}{|c|}{$\mu(x,y)=\cos(x)\cos(4y)(1+\exp(x+y))$ - GMRES} \\
\hline
     &  \multicolumn {2}{|c|}{$\alpha=1.2$} &  \multicolumn {2}{|c|}{$\alpha=1.4$} &  \multicolumn {2}{|c|}{$\alpha=1.6$} &  \multicolumn {2}{|c|}{$\alpha=1.8$}\\
$n$ & - & $P_{\mi{n}}$ & - & $P_{\mi{n}}$ & - & $P_{\mi{n}}$ & - & $P_{\mi{n}}$   \\
   \hline
$2^4$ & 40 & 15 & 45 & 13 & 52 & 11 & 59 & 10 \\
$2^5$ & 62 & 15 & 73 & 13 & 91 & 11 & 111 & 10 \\
$2^6$ & 94 & 16 & 120 & 13 & 158 & 12 & 205 & 10 \\
$2^7$ & 143 & 16 & 196 & 14 & 274 & 12 & 381 & 10 \\
$2^8$ & 216 & 16 & 319 & 14 & 476 & 12 & $>$500 & 10 \\
$2^9$ & 328 & 16 & $>$500 & 14 & $>$500 & 12 & $>$500 & 10 \\
   \hline
  \end{tabular}
  \caption{Number of preconditioned PCG/GMRES iterations to solve the linear system for increasing dimension $n^2$  till $tol =1.e-11$ - $\Omega=[-1,1]^2$.} \label{tab:pcggmres_flag10}
\end{table}
\begin{table}
  \centering
  \footnotesize
  \begin{tabular}{|c |c|c |c|c |c|c |c|c|}
\hline
\multicolumn {9}{|c|}{$\mu(x,y)=\exp(\hat{i}(x+4y))$} \\
\hline
     &  \multicolumn {2}{|c|}{$\alpha=1.2$} &  \multicolumn {2}{|c|}{$\alpha=1.4$} &  \multicolumn {2}{|c|}{$\alpha=1.6$} &  \multicolumn {2}{|c|}{$\alpha=1.8$}\\
$n$ & - & $P_{\mi{n}}$ & - & $P_{\mi{n}}$ & - & $P_{\mi{n}}$ & - & $P_{\mi{n}}$   \\
   \hline
$2^4$ & 39 & 11 & 45 & 10 & 52 & 9 & 60 & 8 \\
$2^5$ & 60 & 12 & 74 & 11 & 90 & 10 & 110 & 8 \\
$2^6$ & 91 & 12 & 120 & 11 & 157 & 10 & 204 & 8 \\
$2^7$ & 138 & 12 & 195 & 11 & 272 & 10 & 380 & 9 \\
$2^8$ & 209 & 13 & 317 & 12 & 476 & 10 & $>$500 & 9 \\
$2^9$ & 319 & 13 & $>$500 & 12 & $>$500 & 11 & $>$500 & 9 \\
   \hline
  \end{tabular}
  \caption{Number of preconditioned GMRES iterations to solve the linear system for increasing dimension $n^2$  till $tol =1.e-11$ - $\Omega=[-1,1]^2$.} \label{tab:gmres_flag11}
\end{table}
\begin{table}
  \centering
  \footnotesize
  \begin{tabular}{|c |c|c |c|c |c|c |c|c|}
\hline
\multicolumn {9}{|c|}{$\mu(x,y)=-2+\exp(\hat{i}(3x+2y))$} \\[2pt]
\hline
     &  \multicolumn {2}{|c|}{$\alpha=1.2$} &  \multicolumn {2}{|c|}{$\alpha=1.4$} &  \multicolumn {2}{|c|}{$\alpha=1.6$} &  \multicolumn {2}{|c|}{$\alpha=1.8$}\\
$n$ & - & $P_{\mi{n}}$ & - & $P_{\mi{n}}$ & - & $P_{\mi{n}}$ & - & $P_{\mi{n}}$   \\
   \hline
$2^4$ & 47 & 15 & 51 & 13 & 57 & 11 & 63 & 10 \\
$2^5$ & 73 & 16 & 84 & 13 & 99 & 12 & 117 & 10 \\
$2^6$ & 111 & 16 & 137 & 14 & 171 & 12 & 216 & 10 \\
$2^7$ & 168 & 16 & 222 & 14 & 298 & 12 & 402 & 10 \\
$2^8$ & 256 & 16 & 360 & 14 & $>$500 & 12 & $>$500 & 10 \\
$2^9$ & 387 & 17 & $>$500 & 14 & $>$500 & 12 & $>$500 & 10 \\
   \hline
  \end{tabular}
  \caption{Number of preconditioned GMRES iterations to solve the linear system for increasing dimension $n^2$  till $tol =1.e-11$ - $\Omega=[-1,1]^2$.} \label{tab:gmres_flag12}
\end{table}


\section{Conclusions}\label{sec:final}

In this work we considered a Helmholtz equation with fractional Laplace operator approximated by ad hoc centered differences with variable wave number $\mu(x,y)$, in great generality including the case in which the function is complex-valued. Several results concerning distribution and clustering have been derived, both in the sense of the eigenvalues and singular values. A wide set of numerical experiments and visualizations has been reported and discussed. The results strongly agree with the theoretical forecasts presented in the theoretical part.

Many more intricate cases can be treated using the same type of theoretical apparatus, including the GLT theory, tools like the approximating class of sequences and that of non-Hermitian perturbations of Hermitian matrix-sequences. In particular, future work can be focused on the following topics
\begin{itemize}
\item A $p$ order finite element approximation in space dimension $d$ leads to $d$-level Toeplitz matrices having a matrix-valued generating function, which is also the GLT symbol. In this case the parameter $r$ takes the form $r=p^d$ \cite{GSS}. 
\item Along the same lines of the previous item, we could consider discontinuous Galerkin techniques of degree $p$ with  $r=(p+1)^d$ \cite{GLT-blockdD}, or isogeometric analysis with polynomial degree $p$ and regularity $k\le p-1$ with  $r=(p-k)^d$ \cite{GLT-blockdD,bene1,bene2,GLT-exposition-eng}.
\item The preceding two items could be considered also in the case where the main operator, i.e. the fractional Laplacian, is defined on a non-Cartesian $d$-dimensional domain, or equipped with variable coefficients, or with approximations on graded grids. We observe that cases considered in the current item lead to $r\times r$ matrix-valued symbols in $t=2d$ independent variables, as in the general GLT theory \cite{Barb,GLT-blockdD,GLT-BookII}. In fact, the related GLT theory is already available \cite{Barb,GLT-block1D,GLT-blockdD,BS,GLT-BookI,GLT-BookII}, while the other tools do not depend from a specific structure of the involved matrix sequences. Of course when a vector PDE or a vector FDE is considered, the number $r$ has to be multiplied by the number of scalar differential equations (see e.g. \cite{doro,axio} and references therein).
\end{itemize}

\section*{Acknowledgements}

The work of Andrea Adriani, Rosita L. Sormani, Stefano Serra-Capizzano, Cristina Tablino-Possio is supported by GNCS-INdAM.
The work of Rolf Krause and Stefano Serra-Capizzano is funded from the European High-Performance Computing Joint Undertaking  (JU) under grant agreement No 955701. The JU receives support from the European Union’s Horizon 2020 research and innovation programme and Belgium, France, Germany, Switzerland.
Furthermore, Stefano Serra-Capizzano is grateful for the support of the Laboratory of Theory, Economics and Systems – Department of Computer Science at Athens University of Economics and Business.

\section*{Data availability statement}

The authors confirm that the data supporting the findings of this study are available within the article.




\begin{thebibliography}{MMMMM}

\bibitem{AST}
A. Adriani, S. Serra-Capizzano, C. Tablino-Possio, Clustering analysis and preconditioned Krylov solvers for the approximated Helmholtz equation and fractional laplacian in the case of complex-valued, unbounded variable coefficient wave number $\mu$, submitted.

\bibitem{model-1} C. Bucur, E. Valdinoci, Nonlocal Diffusion and Applications, Springer, Cham, 2016.

\bibitem{BEV}
N. Barakitis, S.-E. Ekström, P. Vassalos, Preconditioners for fractional diffusion equations based on the spectral symbol, Numer. Linear Algebra Appl. 29 (2022), no. 5, e2441, 22 pp.




\bibitem{Barb}
G. Barbarino,  A systematic approach to reduced GLT, BIT 62 (2022), no. 3, 681–743.

\bibitem{GLT-block1D} G. Barbarino, C. Garoni, S. Serra-Capizzano, {Block generalized locally Toeplitz sequences: theory and applications in the unidimensional case}, Electr. Trans. Numer. Anal. 53 (2020), 28-112.

\bibitem{GLT-blockdD} G. Barbarino, C. Garoni, S. Serra-Capizzano, {Block generalized locally Toeplitz sequences: theory and applications in the multidimensional case},  Electr. Trans. Numer. Anal. 53 (2020), 113--216.

\bibitem{non-herm-perturb} G. Barbarino, S. Serra-Capizzano,  Non-Hermitian perturbations of Hermitian matrix-sequences and applications to the spectral analysis of the numerical approximation of partial differential equations, Numer. Linear Algebra Appl. 27 (2020), no. 3, e2286, 31 pp.


\bibitem{BS}
B. Beckermann, S. Serra-Capizzano, On the asymptotic spectrum of finite element matrix sequences, SIAM J. Numer. Anal. 45 (2007), no. 2, 746–769.

\bibitem{bene1}
P. Benedusi, P. Ferrari, C. Garoni, R. Krause, S. Serra-Capizzano, Fast parallel solver for the space-time IgA-DG discretization of the diffusion equation, J. Sci. Comput. 89 (2021), no. 1, Paper No. 20, 21 pp.

\bibitem{bene2}
P. Benedusi, C. Garoni, R. Krause, X. Li, S. Serra-Capizzano, Space-time FE-DG discretization of the anisotropic diffusion equation in any dimension: the spectral symbol, SIAM J. Matrix Anal. Appl. 39 (2018), no. 3, 1383–1420.


\bibitem{Bhatia-book}
R.~Bhatia, { {M}atrix {A}nalysis}, { Graduate Texts in Mathematics}, 169, Springer-Verlag, New York, 1997.


\bibitem{BC}
D. Bini, M. Capovani, Spectral and computational properties of band symmetric Toeplitz matrices, Linear Algebra Appl. 52/53 (1983), 99–126.


\bibitem{BG-extr}
A. Böttcher, S. Grudsky, On the condition numbers of large semi-definite Toeplitz matrices, Linear Algebra Appl. 279 (1998), no. 1-3, 285–301.

\bibitem{CN} R. Chan, M. Ng,  Conjugate gradient methods for Toeplitz systems, SIAM Rev. 38 (1996), no. 3, 427–482.

\bibitem{17} M. D’Elia, M. Gunzburger, The fractional Laplacian operator on bounded domains as a special case of the nonlocal diffusion operator, Comput. Math. Appl. 66 (2013), 1245–1260.

\bibitem{DB}
F. Di Benedetto, Preconditioning of block Toeplitz matrices by sine transforms, SIAM J. Sci. Comput. 18 (1997), no. 2, 499–515.

\bibitem{DiB-S}
F. Di Benedetto, S.Serra-Capizzano, Optimal multilevel matrix algebra operators, Linear and Multilinear Algebra 48 (2000), no. 1, 35–66.




\bibitem{doro}
A. Dorostkar, M. Neytcheva, S. Serra-Capizzano, Spectral analysis of coupled PDEs and of their Schur complements via generalized locally Toeplitz sequences in 2D, Comput. Methods Appl. Mech. Engrg. 309 (2016), 74–105.

\bibitem{13} S.W. Duo, H.W. van Wyk, Y.Z. Zhang, A novel and accurate finite difference method for the fractional Laplacian and the fractional Poisson problem, J. Comput. Phys. 355 (2018), 233–252.

\bibitem{FaTi}
D. Fasino, P. Tilli, Spectral clustering properties of block multilevel Hankel matrices, Linear Algebra Appl. 306 (2000), no. 1-3, 155–163.



\bibitem{19} A. Gabriel, J.P. Borthagaray, A fractional Laplace equation: regularity of solutions and finite element approximations, SIAM J. Numer. Anal. 55 (2017), 472–495.

\bibitem{axio}
C. Garoni, M. Mazza, S. Serra-Capizzano, Block generalized locally Toeplitz sequences: From the theory to the applications,
Axioms 7 (2018), no. 32018, paper 49.


\bibitem{GLT-BookI} C. Garoni, S. Serra-Capizzano, Generalized locally Toeplitz sequences: theory and applications. Vol. I, Springer, Cham, 2017.
\bibitem{GLT-BookII} C. Garoni,  S. Serra-Capizzano, Generalized locally Toeplitz sequences: theory and applications. Vol. II, Springer, Cham, 2018.

\bibitem{GSS} C. Garoni, S. Serra-Capizzano, D. Sesana, pectral analysis and spectral symbol of $d$-variate ${\mathbb Q}_p$ Lagrangian FEM stiffness matrices, SIAM J. Matrix Anal. Appl. 36 (2015), no. 3, 1100–1128.

\bibitem{GLT-exposition-eng} C. Garoni, H. Speleers, S.-E. Ekstr\"om, S. Serra-Capizzano, T.J.R. Hughes, Symbol-based analysis of finite element and isogeometric B-spline discretizations of eigenvalue problems: exposition and review, Arch. Comput. Methods Eng. 26 (2019), no. 5, 1639–1690.

\bibitem{9} P. Gatto, J.S. Hesthaven, Numerical approximation of the fractional Laplacian via hp-finite elements, with an application to image denoising, J. Sci. Comput. 65 (2015), 249–270.


\bibitem{18} C. Glusa, H. Antil, M. D’Elia, B. van Bloemen Waanders, C.J. Weiss, A Fast Solver for the Fractional Helmholtz Equation, SIAM J. Sci. Comput. 43 (2021), A1362–A1388.

\bibitem{orth-distribution2}
L. Golinskii, S. Serra-Capizzano, The asymptotic properties of the spectrum of nonsymmetrically perturbed Jacobi matrix sequences, J. Approx. Theory 144 (2007), no. 1, 84–102.


\bibitem{16} Z.P. Hao, Z.Q. Zhang, Optimal regularity and error estimates of a spectral Galerkin method for fractional advection-diffusion-reaction equations, SIAM J. Numer. Anal. 58 (2020), 211–233.

\bibitem{20} Z.P. Hao, Z.Q. Zhang, R. Du, Fractional centered difference scheme for high-dimensional integral fractional Laplacian, J. Comput. Phys. 424 (2021), 109851.



\bibitem{14} Y.H. Huang, A. Oberman, Numerical methods for the fractional Laplacian: A finite difference quadrature approach, SIAM J. Numer. Anal. 52 (2014), 3056–3084.

\bibitem{KO}
T. Kailath, V. Olshevsky, Displacement structure approach to discrete-trigonometric-transform based preconditioners of G. Strang type and of T. Chan type, SIAM J. Matrix Anal. Appl. 26 (2005), no. 3, 706–734.


\bibitem{11} M. Kwasnicki, Ten equivalent definitions of the fractional Laplace operator, Fract. Calc. Appl. Anal. 20 (2017), 7–51.

\bibitem{orth-distribution1}
A. Kuijlaars, S. Serra-Capizzano, Asymptotic zero distribution of orthogonal polynomials with discontinuously varying recurrence coefficients, J. Approx. Theory 113 (2001), no. 1, 142–155.


\bibitem{10} N. Laskin, Fractional quantum mechanics and Levy path integrals, Phys. Lett. A 268 (2000), 298–305.

\bibitem{tau-prec} T.-Y. Li,  F. Chen,  H.W. Sun, T. Sun,  Preconditioning Technique Based on Sine Transformation for Nonlocal Helmholtz Equations with Fractional Laplacian, J. Sci. Comput. 97 (2023), no. 1, 17.

\bibitem{12} A. Lischke, G. Pang, M. Gulian, F. Song, C. Glusa, X. Zheng, Z. Mao, W. Cai, M.M. Meerschaert, M. Ainsworth, G.E. Karniadakis, What is the fractional Laplacian? A comparative review with new results, J. Comput. Phys. 404 (2020), 109009.

\bibitem{model-2} M. Meerschaert, A. Sikorskii, Stochastic Models for Fractional Calculus, De Gruyter, 2019.

\bibitem{15} V. Minden, L.X. Ying, A simple solver for the fractional Laplacian in multiple dimensions, SIAM J. Sci. Comput. 42 (2020), A878–A900.



\bibitem{tau-essential} D. Noutsos, S. Serra-Capizzano, P. Vassalos,  Essential spectral equivalence via multiple step preconditioning and applications to ill conditioned Toeplitz matrices, Linear Algebra Appl. 491 (2016), 276–291.




\bibitem{model-3} I. Podlubny, Fractional Differential Equations, Academic Press, New York, 1999.


\bibitem{S-extr2} S. Serra-Capizzano, On the extreme spectral properties of Toeplitz matrices generated by $L^1$ functions with several minima/maxima, BIT 36 (1996), no. 1, 135–142.
\bibitem{S-extr1} S. Serra-Capizzano, On the extreme eigenvalues of Hermitian (block) Toeplitz matrices, Linear Algebra Appl. 270 (1998), 109–129.




\bibitem{tau-theory1} S. Serra-Capizzano, Superlinear PCG methods for symmetric Toeplitz systems, Math. Comp. 68 (1999), no. 226, 793–803.
\bibitem{tau-theory2} S. Serra-Capizzano, Toeplitz preconditioners constructed from linear approximation processes. SIAM J. Matrix Anal. Appl. 20 (1999), no. 2, 446–465.

\bibitem{taud2}
S. Serra-Capizzano, Spectral behavior of matrix sequences and discretized boundary value problems, Linear Algebra Appl. 337 (2001), 37–78.

\bibitem{nega-gen}
S. Serra-Capizzano, Matrix algebra preconditioners for multilevel Toeplitz matrices are not superlinear. Special issue on structured and infinite systems of linear equations. Linear Algebra Appl. 343/344 (2002), 303–319. 

\bibitem{prec-neg} S. Serra Capizzano, E. Tyrtyshnikov,  How to prove that a preconditioner cannot be superlinear, Math. Comp. 72 (2003), no. 243, 1305–1316.


\bibitem{Van} C. Van Loan, Computational frameworks for the fast Fourier transform, Frontiers in Applied Mathematics, 10. Society for Industrial and Applied Mathematics (SIAM), Philadelphia, PA, 1992.

\bibitem{model-4} J.L. Vazquez, The mathematical theories of diffusion: Nonlinear and fractional diffusion, in: Nonlocal and Nonlinear Diffusions and Interactions: New Methods and Directions, Springer, 2017, 205–278.







\end{thebibliography}
\end{document}